\documentclass{article}

\usepackage{amsmath}
\usepackage{amssymb}
\usepackage{graphicx}
\usepackage{geometry} 
\usepackage{pgf,tikz,pgfplots}
\usepackage{mathrsfs}
\usetikzlibrary{arrows}
\usepackage{amsthm}
\usepackage[utf8]{inputenc}
\usepackage{titlesec}
\usepackage{multicol}
\usepackage[normalem]{ulem}
\usepackage[T1]{fontenc} 
\usepackage[utf8]{inputenc} 
\usepackage{array}
\usepackage{calc}
\usepackage{tikz-3dplot}
\usepackage{authblk} 

\definecolor{rouge}{RGB}{204, 15, 39}
\definecolor{bleu}{RGB}{60, 86, 164}
\definecolor{vert}{RGB}{10, 195, 118}
\definecolor{orange}{RGB}{245, 85, 63}
\definecolor{gris}{RGB}{69, 94, 93}
\definecolor{noir}{RGB}{0, 0, 0}
\definecolor{jaune}{RGB}{255, 195, 42}

\begingroup 
\makeatletter 
\@for\theoremstyle:=definition,remark,plain\do{%
	\expandafter\g@addto@macro\csname th@\theoremstyle\endcsname{%
		\addtolength\thm@preskip\parskip 
	}%
} 
\endgroup 

\newenvironment{itemize*}{\begin{itemize}\setlength{\itemsep}{-2pt} \vspace{-0.17cm}\end{itemize}}

\newtheorem*{thm*}{Theorem}
\newtheorem*{pties*}{Properties}
\newtheorem*{prop*}{Proposition}
\newtheorem*{defn*}{Definition}
\newtheorem*{defnthm*}{Definition and Theorem}

\newtheorem{thm}{Theorem}[section]
\newtheorem{lem}[thm]{Lemma}

\newtheorem{prop}[thm]{Proposition}
\newtheorem{cor}[thm]{Corollary}
\newtheorem{defn}[thm]{Definition}
\newtheorem{defnthm}[thm]{Definition and Theorem}

\theoremstyle{remark}
\newtheorem{rem}[thm]{Remark}
\newtheorem{ex2}[thm]{Example}
\newenvironment{ex}{\begin{ex2}}{\hfill $\diamondsuit$ \end{ex2}}

\newtheorem{fex2}[thm]{Fundamental example}
\newtheorem{cex2}[thm]{Counter example}
\newtheorem{iex2}[thm]{Interesting examples}

\newtheorem*{pro}{\textit{Proof}}
\newtheorem*{proth}{\textit{Proof of the theorem}}
\newtheorem*{propr}{\textit{Proof of the proposition}}
\newtheorem*{proPr}{\textit{Proof of Proposition}}
\newtheorem{nota}[thm]{Notation}
\newenvironment{prr}{\begin{pro}}{\hfill $\square$ \end{pro}}

\newenvironment{prPr}{\begin{proPr}}{\hfill $\square$ \end{proPr}}

\makeatletter
\newtheoremstyle{indented}
{4pt}
{6pt}
{\addtolength{\@totalleftmargin}{33.3pt}
	\addtolength{\linewidth}{-33.3pt}
	\parshape 1 33.3pt \linewidth }
{-33.3pt}
{}
{.}
{.5em}
{}
\makeatother
\theoremstyle{indented}
\newtheorem*{pr2}{\textit{Proof}}
\newtheorem*{proft2}{\textit{Proof of the theorem}}
\newtheorem*{profp2}{\textit{Proof of the proposition}}

\makeatletter
\newtheoremstyle{indentedlem}
{4pt}
{6pt}
{\itshape\addtolength{\@totalleftmargin}{33.3pt}
	\addtolength{\linewidth}{0pt}
	\parshape 1 33.3pt \linewidth }
{}
{}
{.}
{.3em}
{}
\makeatother
\theoremstyle{indentedlem}

\makeatletter
\newtheoremstyle{indented2}
{2pt}
{3pt}
{\addtolength{\@totalleftmargin}{0pt}
	\addtolength{\linewidth}{-33.3pt}
	\parshape 1 66.6pt \linewidth }
{-33.3pt}
{}
{.}
{.5em}
{}
\makeatother
\theoremstyle{indented2}
\newtheorem*{prl2}{\textit{Proof}}

\newcommand{\bC}{\mathbb{C}}

\newcommand{\bF}{\mathbb{F}}

\newcommand{\bN}{\mathbb{N}}

\newcommand{\bP}{\mathbb{P}}
\newcommand{\bQ}{\mathbb{Q}}
\newcommand{\bR}{\mathbb{R}}
\newcommand{\bS}{\mathbb{S}}

\newcommand{\bZ}{\mathbb{Z}}

\newcommand{\fM}{\mathfrak{M}}

\newcommand{\fS}{\mathfrak{S}}

\newcommand{\fX}{\mathfrak{X}}

\newcommand{\sB}{\mathscr{B}}

\newcommand{\sH}{\mathscr{H}}

\newcommand{\sO}{\mathscr{O}}

\newcommand{\su}{\textup{sup}}

\newcommand{\vo}{\textup{vol}}

\newcommand{\bi}{\textup{Big}}
\newcommand{\rk}{\textup{rk}}

\newcommand{\Proj}{\textup{Proj}}

\newcommand{\Nef}{\textup{Nef}}
\newcommand{\Bigu}{\textup{Big}}
\newcommand{\Eff}{\textup{Eff}}

\newcommand{\mo}{\textup{Mov}}
\newcommand{\Sym}{\textup{Sym}}
\newcommand{\cc}{\textup{cc}}
\newcommand{\fsl}{\textup{fs}}

\usepackage{setspace}
\definecolor{link}{RGB}{71,81,114}
\usepackage{hyperref}                 
\hypersetup{
	colorlinks = true,
	allcolors=link}

\title{Newton--Okounkov bodies of curve classes}
\author{Lucie Devey}
\date{}
\newcommand{\Addresses}{{
		\bigskip
		\footnotesize\textsc{Institut für Mathematik, Goethe–Universität Frankfurt, 60325 Frankfurt am Main, Germany\\
		Institut Fourier, Université Grenoble Alpes, 38400 Saint Martin d’Hères, France}\par\nopagebreak
		\textit{E-mail address}, \texttt{lucie.devey@univ-grenoble-alpes.fr}
}}

\begin{document}

	\maketitle

	\begin{abstract}
		The purpose of the paper is to initiate the development of the theory of Newton--Okounkov bodies of curve classes. 
		
		Our definition is based on making a fundamental property of Newton--Okounkov bodies hold also in the curve case: the volume of the Newton--Okounkov body of a curve is a volume-type function of the original curve.
		This construction allows us to conjecture a new relation between Newton--Okounkov bodies, we prove it in certain cases.
	\end{abstract}

\section*{Introduction}

In the late 1990s Okounkov associated a convex body $\Delta_{Y_\bullet}(D)$ to any ample divisor $D$ on a projective variety $X$ depending on the choice of an admissible flag $Y_\bullet$ and used its geometric properties to explore the sections $H^0(X,\sO_X(kD))$ for large values of $k$. In 2008, Lazarsfeld--Musta\c{t}\u{a} \cite{ref4} and Kaveh--Khovanskii \cite{ref7} simultaneously extended this contruction to any big divisor and developed the theory of Newton--Okounkov bodies. 

This object is a vast generalisation of the Newton polytope of a hypersurface and has come to be known as the 'Newton--Okounkov body' of $D$. This article is an attempt to generalise the theory defining Newton--Okounkov bodies of curve classes instead of divisor classes and to discover new relations between Newton--Okounkov bodies. 

The interest in Newton--Okounkov bodies lies in the fact that they encode information on the asymptotic behaviour of $H^0(X,mD)$ as $m\to\infty$ and have been proved to be very efficient in providing simple geometric proofs of difficult results such as Fujita's approximation theorem (\cite{ref4} Theorem 3.5) or the log-concavity of the volume of a divisor (\cite{ref3} Theorem 1.6). 

The Newton--Okounkov body of a divisor can also reveal information about its Seshadri constants, invariants introduced by Demailly and measuring the positivity of the divisor. In \cite{ref25}, Ito furnished a lower bound on Seshadri constants and Dumnicki--Küronya--Maclean--Szemberg exhibited in \cite{ref24} an unexpected relation between SHGH conjecture and rationality problems for Seshadri constants using Newton--Okounkov bodies. 

As highlighted in the survey article \cite{ref27}, there is a strong relationship
between positivity of divisors and the geometry of Newton–Okounkov bodies. This connection persists locally: \cite{ref28} and \cite{ref29} relate local positivity of line bundles, for instance jet-separation, to Newton–Okounkov bodies attached to infinitesimal
flags.

Even though Newton--Okounkov theory for divisors turned out to be very fruitful, there has so far been no visible attempt to construct higher-codimensional generalisations.

One of the most remarkable features of Newton--Okounkov bodies is that the volume of a divisor
$$ \vo(D)=\textup{lim\,sup}_{m\to\infty}\frac{h^0(X,mD)}{m^n/n!}$$
is exactly the euclidean volume of its Newton--Okounov body for any choice of admissible flag.
We want to propose a definition for Newton--Okounkov bodies of curve classes which still have this concrete and geometrical property. 
Lehmann and Xiao (\cite{ref3}) constructed a Legendre--Fenchel type transform $\fM$ on curve classes by taking the dual of the volume function :
$$\fM(\alpha)=\inf_{\substack{A\text{ big and movable}\\\text{divisor class}}}\bigg(\frac{A\cdot\alpha}{\vo(A)^\frac1n}\bigg)^\frac n{n-1}\ .$$
Our definition of Newton--Okounkov bodies for curve classes will satisfy $$\vo_{\bR^n}\left(\Delta(\alpha)\right)=\fM(\alpha)\ .$$

If $X$ is a surface then curves are divisors and consequently Newton--Okounkov bodies of curves are already defined. Moreover an explicit description of Newton--Okounkov bodies can be given based on the Zariski decomposition for divisors (see \cite{ref4} Theorem 6.4 or \cite{ref15} Section 2 for more details). 

Lehmann and Xiao have defined a Zariski-type decomposition for curve classes which we will call the 'movable Zariski decomposition', based on the volume-type function $\fM$ (see \cite{ref3} Theorem 3.12).
Assume that the curve class $\alpha$ is movable (or equivalently by \cite{ref16} the dual of $\overline{\Eff_1}(X)$) and $\fM(\alpha)>0$, Lehmann--Xiao's result states that there exists a unique big and movable divisor class $L$ such that
	$$\alpha=\langle L^{n-1}\rangle$$ 
where $\langle \rangle$ is the positive product (see Boucksom--Favre--Jonsson in \cite{ref23}).
Moreover the infimum arising in the definition of $\fM$ is attained by $L$.

Based on this decomposition, we propose a definition for the Newton--Okounkov body of a curve class $\alpha$ with respect to an admissible flag $Y_\bullet$ on $X$.

\begin{defn*}[\ref{defn}]
	We define the Newton--Okounkov body of a movable curve class $\alpha$ such that $\fM(\alpha)>0$ as
	$$\Delta_{Y_\bullet}(\alpha)=\Delta_{Y_\bullet}(L)\ ,$$
	where $\alpha=\langle L^{n-1}\rangle$ is the movable Zariski decomposition of $\alpha$.
\end{defn*}
We present some formal properties, namely the equality between the volume of a curve class and the volume of its Newton--Okounkov body, and the continuity of Newton--Okounkov bodies of curve classes. 

In the second part of the paper we investigate a potential analogue of the inclusion 	
	$$\Delta_{Y_\bullet}(D_1)+\Delta_{Y_\bullet}(D_2)\subseteq\Delta_{Y_\bullet}(D_1+D_2)\ .$$ 
No such simple inclusion appears to hold for Newton--Okounkov bodies of curves as defined above. Even in the most simple case, when $X=\bP^3$ and the curve classes $\alpha_1=\alpha_2=\langle D^2\rangle$ are given by $D=\sO_X(1)$, the movable Zariski decomposition of the sum $\alpha_1+\alpha_2$ is $\alpha_1+\alpha_2=\langle(\sqrt2D)^2\rangle$ and
$$\Delta(\alpha_1)+\Delta(\alpha_2)=2\Delta(D)\nsubseteq\sqrt2\Delta(D)=\Delta(\alpha_1+\alpha_2)\ .$$ 

A potential replacement for the Minkowski sum in the context of curve classes is the Blaschke sum $\#$ of convex sets (see Definition \ref{bs}). We would like to prove that the inclusion

\begin{equation}\label{incmov}
	\Delta_{Y_\bullet}(\alpha_1)\#\Delta_{Y_\bullet}(\alpha_2)\subseteq\Delta_{Y_\bullet}(\alpha_1+\alpha_2) \qquad\qquad \tag{$*_{\textup{inc}}$}
\end{equation}
holds. In this case we would have proved that, for any movable divisor $L_1$ and $L_2$ on $X$, denoting by $L_3$ the unique movable divisor satisfying $\langle L_3^{n-1}\rangle=\langle L_1^{n-1}\rangle+\langle L_2^{n-1}\rangle$,
\begin{equation*}
	\Delta_{Y_\bullet}(L_1)\#\Delta_{Y_\bullet}(L_2)\subseteq\Delta_{Y_\bullet}(L_3)\ .
\end{equation*}
Since the definition of the Blaschke sum is entirely based on the area measure of convex sets (see Definition \ref{aream}), this inclusion would be one of the first result we have on the boundary of Newton--Okounkov bodies. 

One motivation for considering Blaschke sums is that the sum of curve classes and the Blaschke sum of convex bodies satisfy the volume formula (see \cite{ref10} Section 7.A) 
\begin{align*}
	&\vo(K\#L)^{\frac{n-1}n}\geq\vo(K)^{\frac{n-1}n}+\vo(L)^{\frac{n-1}n}\\
	\text{and}\quad &\fM(\alpha_1+\alpha_2)^{\frac{n-1}n}\geq\fM(\alpha_1)^{\frac{n-1}n}+\fM(\alpha_2)^{\frac{n-1}n}\ .
\end{align*}

In this paper we prove (\ref{incmov}) in the case of surfaces, of homothetic curve classes and of toric varieties. It turns out that in the two last cases (\ref{incmov}) is an equality. The following result is really powerful but would require the Blaschke sum to be continuous.
\begin{prop*}[\ref{reduction}]
	We assume the continuity of the Blaschke sum and consider Newton--Okounkov bodies with respect to any maximal rank valuation (not only flag valuation). 
	
	Then the inclusion (\ref{incmov}) holds for movable curves if and only if it holds for curves of the form $\alpha_i=A_i^{n-1}$ with $A_i$ ample. 
\end{prop*} 
Finally we prove the inclusion (\ref{incmov}) in the following case.
\begin{thm*}[\ref{mthm}]
	Let $X$ be a projective bundle over a curve. Consider two curves $\alpha_1,\alpha_2$ of the form $\alpha_i=A_i^{n-1}$ with $A_i$ ample. Then we have that
	$$\Delta(\alpha_1)\#\Delta(\alpha_2)=\Delta(\alpha_1+\alpha_2)\ .$$
\end{thm*}

In the \ref{defin}$^{\text{st}}$ section we propose the definition of the Newton--Okounkov body of a curve class and present some of its properties. 

In Section \ref{mot}, we discuss the Blaschke sum and present various analogies with Newton--Okounkov bodies. 

In Section \ref{cont} we prove Proposition \ref{reduction}.

Section \ref{nob,cones} describes Newton--Okounkov bodies of curve classes when $X$ is a projective bundle over curves, recalls the description of the cones of divisor and curve classes and explains how the Newton--Okounkov body allows to visualize the positivity. 

In section \ref{Inc}, we prove Theorem \ref{mthm}.

\subsection*{Acknowledgements}
The author thanks her supervisors Catriona Maclean and Alex Küronya for their many corrections and enlightening discussions. The author is grateful to Stefano Urbinati for his helpful comments.

This work has been partially supported by the LabEx PERSYVAL-Lab (ANR-11-LABX-0025-01) funded by the French program Investissement d’avenir. This work has been partially supported by the Deutsche Forschungsgemeinschaft (DFG, German Research Foundation) TRR 326 Geometry and Arithmetic of Uniformized Structures, project number 444845124.

\section{Definition of Newton--Okounkov bodies of curves}\label{defin}

We start with some general definitions. The theory of Newton--Okounkov bodies was developed simultaneously by Lazarsfeld--Musta\c{t}\u{a} \cite{ref4} and Kaveh--Khovanskii \cite{ref7}. We propose an extension of this definition to curve classes.

Throughout this chapter, $X$ will be a projective variety of dimension $n$ and $Y_\bullet$ an admissible flag on $X$.

\subsection{Newton--Okounkov bodies of curve classes}

Lehmann and Xiao constructed a Zariski-type decomposition for curve classes generalizing the Zariski decomposition of divisors on surfaces. It will allow us to construct an analogous definition of Newton--Okounkov bodies of curve classes. 

Let us start with the volume function of curves introduced by Lehmann and Xiao (see \cite{ref3}). Recall that the volume function for divisors is given by
$$ \vo(D)=\textup{lim\,sup}_{m\to\infty}\frac{h^0(X,mD)}{m^n/n!}\ .$$ 

\begin{defn}\label{vol}
	The dual volume function $\fM$ on curves is defined by
	$$\fM(\alpha)=\inf_{\substack{A\text{ big and movable}\\\text{divisor class}}}\bigg(\frac{A\cdot\alpha}{\vo(A)^\frac1n}\bigg)^\frac n{n-1}\ .$$
\end{defn}

We now recall the definition of the movable Zariski decomposition from \cite{ref3} Theorem 3.12.

\begin{defnthm} \label{mzd}
	Any movable curve class $\alpha$ with $\fM(\alpha)>0$ is of the form
		$$\alpha=\langle L^{n-1}\rangle$$
	for a unique big and movable divisor class $L$, where $\langle \rangle$ is the positive intersection product (Remark \ref{pip}).
	Moreover the infimum appearing in the definition of $\fM$ is achieved by $L$.
\end{defnthm}

\begin{rem} \label{pip}
	The positive intersection product, introduced by Boucksom--Favre--Jonsson (see \cite{ref23}), is defined for classes on the Riemann-Zariski space $\fX$ of a projective variety $X$ which is the projective limit of all birational models of $X$. A class in $\fX$ is a collection of classes in each birational model of X that are compatible under push-forward. We denote the set of such classes by $N_p(\fX)$. 
	
	A class $L$ is Cartier if and only if there exists a birational model $X_\pi$ of $X$ such that the incarnations of $L$ on higher blow-ups are obtained by pulling-back the incarnation $L_\pi$ of $L$ on $X_\pi$. Such a $\pi$ is called a determination of $L$. 
	
	If $L_1,...,L_p$ are big Cartier divisor classes, then their positive intersection product $\langle L_1,...,L_p\rangle$ is defined as the least upper bound of the set of classes 
	$$(L_1-D_1)\cdot...\cdot(L_p-D_p)\in N^p(\fX)$$ 
	where $D_i$ is an effective Cartier Q-divisor on $X$ such that $L_i-D_i$ is nef.
	
	The most relevant property for us is that if $L_1,L_2,...,L_p\in N^1(\fX)$ are nef Cartier divisor classes then
		$$\langle L_1,...,L_p\rangle=L_1\cdot...\cdot L_p\ .$$
	
	Another essential property is Fujita's theorem (see \cite{ref22}) reformulated in \cite{ref23} Theorem 3.1 into
	$$\forall L\in\bi^1(X),\quad \vo(L) =\langle L^n\rangle\ .$$
\end{rem}

\begin{rem}
	The condition $\fM(\alpha)>0$ is equivalent to having non vanishing intersection with any non-zero movable divisor class (see \cite{ref3} Lemma 3.9).
\end{rem}

\begin{rem}\label{dim2dec}
	We call the decomposition 'movable Zariski' because of the Zariski decomposition of effective divisors (see \cite{ref11})
		$$D=B+\gamma$$ 
	where $B$ is nef, $\gamma$ is a negative cycle and $B\cdot\gamma=0$ from which it follows that the dimension of the linear system $|D|$ is determined by $B$ alone. 
	
	So up to a translation (depending on $\gamma$), the Newton--Okounkov body of $D$ is the Newton--Okounkov body of $B$ (see Paragraph 6.2 \cite{ref4}). 
	
	With the 'movable Zariski decomposition' for curve classes, Lehmann and Xiao extended the decomposition of Zariski to movable curve classes positive with respect to $\fM$.
\end{rem} 
	
This motivates the following definition.

\begin{defn}\label{defn}
	 Consider a movable curve class $\alpha$ on $X$ that satisfies $\fM(\alpha)>0$ and hence has a movable Zariski decomposition $\alpha=\langle L^{n-1}\rangle$. We define the Newton--Okounkov body of $\alpha$ on $X$ with respect to $Y_\bullet$ as
		$$\Delta_{Y_\bullet}(\alpha):=\Delta_{Y_\bullet}(L)\ .$$
\end{defn}

\begin{rem}
	The Newton--Okounkov bodies of curve classes are well defined. Indeed, Newton--Okounkov bodies of big divisors are invariant under numerical equivalence (see Proposition 4.1 \cite{ref4}).
\end{rem}

\begin{ex}
	Take any movable curve class $\alpha$ on $X=\bP^n$ with an admissible flag $Y_\bullet$. The intersection ring of $X$ is $$A(X)=\bZ[H]/(H^n)\ ,$$ so that $\alpha$ can be written as $$\alpha=aH^{n-1}=a\langle H^{n-1}\rangle\ ,$$
	where $a\in\bR$.
	Indeed the cone of nef divisors on $X$ coincides with the cone of movable divisors and taking the positive intersection product of nef divisors corresponds to taking the intersection of these divisors (see Remark \ref{pip}). The Newton--Okounkov body of $\alpha$ is then $$\vspace{-24pt}\Delta_{Y_\bullet}(\alpha)=\Delta_{Y_\bullet}(a^\frac{1}{n-1}H)\ .$$\vspace{6pt}
\end{ex}

\subsection{Formal properties}

Let us fix a movable curve
class $\alpha$ with $\fM(\alpha)>0$, we denote its movable Zariski decomposition by $\alpha=\langle L_\alpha^{n-1}\rangle$. 
All Newton--Okounkov bodies are defined with respect to the fixed flag $Y_\bullet$ so that we omit $Y_\bullet$ in our notation and write $\Delta(\alpha)$. 

We now summarize some important properties of Newton--Okounkov bodies of curves. 

The volume of the Newton--Okounkov body of a divisor computes the volume of the divisor.
The volume of the Newton--Okounkov body of a curve $\alpha$ turns out to be related to a volume of $\alpha$ (defined by Lehmann and Xiao) as well, and is given by a geometric intersection.

\begin{prop}
	Consider a movable curve $\alpha$ such that $\fM(\alpha)>0$. We have the following equalities 
	$$n!\vo\big(\Delta_{Y_\bullet}(\alpha)\big)=\fM(\alpha)=\vo(L_\alpha)=\langle L_\alpha^{n}\rangle\ .$$
\end{prop}

\begin{prr}
	We have that 	
	$$n!\vo(\Delta(\alpha))=\vo(\Delta(L_\alpha))=\vo(L_\alpha)$$ 
	by \cite{ref4} Theorem 2.3. Then Theorem 3.12 of \cite{ref3} states that 
	$$\fM(\alpha)=\vo(L_\alpha)$$
	and the Fujita's theorem of Remark \ref{pip} implies that
	$$\vspace{-24pt}\vo(L_\alpha)=\langle L_\alpha^{n}\rangle\ .$$\vspace{6pt}
\end{prr}

\begin{rem}
	By \cite{ref4} Proposition 4.1, for any divisor $D$ on $X$ and any integer $p>0$, one has 
	$$\Delta(pD)=p\Delta(D)\ .$$ 
	We have an analogous result for any curve class $\alpha$:
	$$\Delta(p\alpha)=p^{\frac1{n-1}}\Delta(\alpha)\ .$$
\end{rem}

Another property is the continuity of Newton--Okounkov bodies of curve classes.

\begin{prop}\label{c0}
	The map $f:\Bigu_1(X)\to\{\text{convex bodies in }\bR^d\}$ defined by 
	$$f:\alpha\mapsto\Delta(\alpha)$$ 
	is continuous with respect to Hausdorff distance $d_H$ on the set of convex subsets in $\bR^d$.
\end{prop}
	
\begin{prr}
	It follows from Theorem 3.15 \cite{ref3} that the class $L_\alpha$ depends continuously on $\alpha$, therefore it is enough to show that the map $g:\Bigu^1(X)\to\{\text{convex bodies in }\bR^d\}$ defined by 
		$$g:D\mapsto\Delta(D)$$ 
	is continuous with respect to the Hausdorff distance $d_H$ on $\bR^d$. This notion of continuity for Newton--Okounkov bodies follows from the convexity of the global Newton--Okounkov body (see \cite{ref4} Theorem 4.5).
	
	To any point $x$ in the big cone we may associate an arbitrary point $\phi(x)$ contained in the Newton--Okounkov body $\Delta(x)$. We may therefore define a map $$\phi:NS(X)\to\bR^d\ .$$
	
	Now, consider a big divisor $x$ in and a ball $B$ of radius $r$ around $x$ contained in the big cone. Take a point $y$ in $B$ at distance $\delta$ from $x$ and set $w=x-\left(\frac{y-x}\delta\right)r$.
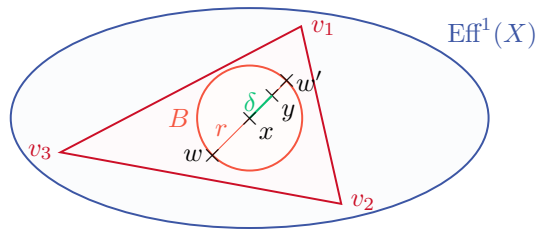
	\begin{center}\begin{tikzpicture}[scale=.75]
		\draw [color=bleu, line width=1.5pt] (0,0) ellipse (120pt and 55pt);
		\fill [color=bleu!2] (0,0) ellipse (120pt and 55pt);
		\draw [color=rouge, line width=1.5pt] (-3.3,-.6) node [left] {$v_3$} -- (1.6,-1.5) node [right] {$v_2$} -- (.9,1.6) node [right] {$v_1$} -- cycle;
		\fill [color=rouge!2] (-3.3,-.6) -- (1.6,-1.5) -- (.9,1.6) -- cycle;
		\draw [color=orange, line width=1.5pt] (0,0) circle (26pt);
		\fill [color=orange!2] (0,0) circle (26pt);
		\draw [color=orange] (-.66,-.66) -- (.66,.66);
		\draw (-.66,-.66) node [left] {$w$};
		\draw (.66,.66) node [right] {$w'$};
		\draw (-.66,-.66) node {$\times$};
		\draw (.66,.66) node {$\times$};
		\draw (.4,.4) node {$\times$};
		\draw [color=vert, line width=.8pt] (0,0) -- (.4,.4);
		\draw (.4,.4) node [below right] {$y$};
		\draw [color=vert] (.29,.29) node [left] {$\delta$};
		\draw (0,0) node {$\times$} node [below right] {$x$};
		\draw [color=orange] (-.22,-.22) node [left] {$r$};
		\draw [color=bleu] (4.3,1.5) node {$\Eff^1(X)$};
		\draw [color=orange] (-.9,0) node [left] {$B$};
	\end{tikzpicture}\end{center}
	Let $q$ be any point in $\Delta(y)$. By convexity of the global Newton--Okounkov body, we have that
		$$s:=\frac{r}{r+\delta}q+\frac{\delta}{r+\delta}\phi(w)\in\Delta(x)\ .$$
	Moreover, if $v_1$, $v_2$ and $v_3$ are big divisor classes forming a triangle containing $B$, by convexity of the global Newton--Okounkov body, every Newton--Okounkov body in the triangle is bounded by $M$ and
		$$d(s,q)\leq\frac \delta{r+\delta}|q|+\frac{\delta}{r+\delta}|\phi(w)|\leq \frac{2M\delta}r\:\underset{\delta\to0}\longrightarrow0$$
	where $M=\su\left\{\textup{bound}(\Delta(v_1)),\textup{bound}(\Delta(v_2)),\textup{bound}(\Delta(v_3))\right\}$.
	
	On the other hand, take $q'$ is a point in $\Delta(x)$ and set $w'=x+(\frac{y-x}\delta)r$. Then the point
		$$s':=\frac{q'(r-\delta)}r+\frac{\phi(w')\delta}r$$
	is contained in $\Delta_{Y_\bullet}(y)$ and 
		$$d(s',q')\leq\frac{\delta}r|q'|+\frac{\delta}r|\phi(w')|\leq\frac{2M\delta} r\:\underset{\delta\to0}\longrightarrow0\ .$$
	
	In conclusion we proved that 	
		$$d_H\big(\Delta(x),\Delta(y)\big)=\textup{max}\left\{\su_{q\in\Delta(x)}d(q,\Delta(y)),\su_{q\in\Delta(y)}d(\Delta(x),q)\right\}\leq2M\frac\delta r$$
	where $M$ depends only on $r$ so that $d_H\big(\Delta(x),\Delta(y)\big)$ tends to $0$ with $\delta$.
\end{prr}

\section{Blaschke sum}\label{mot}

Our long term goal is to study curve classes on a projective complex variety $X$ by means of Newton--Okounkov theory. 
In particular, we would like to study the links between the Newton--Okounkov body $\Delta(\alpha_1+\alpha_2)$ of a sum of curve classes $\alpha_1,\alpha_2$ and the Newton--Okounkov bodies $\Delta(\alpha_1),\Delta(\alpha_2)$ of $\alpha_1$ and $\alpha_2$. 

For divisors, multiplication of sections leads to the inclusion \begin{equation}\tag{$*_{\text{incdiv}}$}\Delta(D_1)+\Delta(D_2)\subseteq\Delta(D_1+D_2)\ . \end{equation}
We will study the following question: can the operation $\#$ (see Definition \ref{bs}) replace the Minkowski sum $+$ in the inclusion ($*_{\text{incdiv}}$) to obtain 
	\begin{equation}\tag{\ref{incmov}}
		\Delta(\alpha_1)\#\Delta(\alpha_2)\subseteq\Delta(\alpha_1+\alpha_2)\ .
	\end{equation}

The Blaschke sum $\#$ constructs from two convex bodies a third one : in the special case where we consider the Blaschke sum of two polytopes $P_1$ and $P_2$, the output $P_3$ is a polytope such that every face of $P_3$ has volume equal to the sum of the volumes of the parallel faces in $P_1$ and $P_2$. 

In fact, while Minkowski sum adds the volumes of dimensional $1$-faces, Blaschke sum adds the volume of codimensional $1$-faces. So we can see Blaschke sum as a dual operation to Minkowski sum. (A formal definition can be found in the next subsection \ref{defnbs}.) 

We are motivated to study the Blaschke sum in this context with a view towards relating $\Delta(\alpha_1)$, $\Delta(\alpha_2)$ and $\Delta(\alpha_1+\alpha_2)$ by the inequalities 
\begin{align*}
	&\vo(K)^{\frac{n-1}n}+\vo(L)^{\frac{n-1}n}\leq\vo(K\#L)^{\frac{n-1}n}\\
	\text{and}\quad &\vo\left(\Delta(\alpha_1)\right)^{\frac{n-1}n}+\vo\left(\Delta(\alpha_2)\right)^{\frac{n-1}n}\leq\vo\left(\Delta(\alpha_1+\alpha_2)\right)^{\frac{n-1}n}
\end{align*}
 satisfied by both the Blaschke sum of convex bodies $K,L$ and the summation of curve classes $\alpha_1,\alpha_2$ (see Section 7.A \cite{ref10}).

\subsection{Definition of Blaschke sum}\label{defnbs}

The Blaschke sum sum is defined in terms of the area measure of a convex body.

\begin{defn}\label{aream}
	The area measure $S_{n-1}(K,\cdot)$ of a convex body $K$ is the measure on the sphere $\bS^{n-1}$ defined by $$S_{n-1}(K,U)=\sH^{n-1}\big(g^{-1}(K,U)\big)$$ where $g^{-1}(K,U)$ is the set of points in $\delta K$ at which there is an outer unit normal vector in $U$ and $\sH^k$ is the $k$-dimensional Hausdorff measure in $\bR^n$ (for $k\in\{1,...,n\}$). 
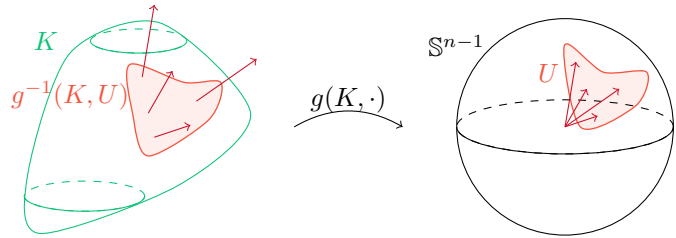
	\begin{center}
		\begin{tikzpicture}[scale=.8]
			\draw [vert] plot [smooth cycle] coordinates {(2.5,1) (1,2.5) (-.5,2) (-1.2,0) (-1,-1) (1.3,-.2)};
			\draw [vert] (-.9,2.2) node {$K$};
			\draw [orange,line width=1] plot [smooth cycle] coordinates {(.8,.3) (.55,1) (.5,1.8) (1.4,1.5) (2,1.5) (1.8,.9)};
			\fill [orange!10] plot [smooth cycle] coordinates {(.8,.3) (.55,1) (.5,1.8) (1.4,1.5) (2,1.5) (1.8,.9)};
			\draw [vert,dashed] (-.26,-.38) ellipse (1cm and 0.25cm);
			\draw [vert] (.74,-.38) arc[x radius=1cm, y radius=.25cm, start angle=0, end angle=-180];
			\draw [vert,dashed] (.64,2.2) ellipse (.8cm and 0.2cm);
			\draw [vert] (1.44,2.2) arc[x radius=.8cm, y radius=.2cm, start angle=0, end angle=-180];
			\draw [rouge,->] (.8,1) -- (1.2,1.7);
			\draw [rouge,->] (.7,1.6) -- (.9,2.8);
			\draw [rouge,->] (1.6,1.2) -- (2.6,1.9);
			\draw [rouge,->] (.9,.6) -- (1.5,.8);
			\draw [orange] (-.5,1.3) node {$g^{-1}(K,U)$};
		\end{tikzpicture}
		\quad
		\begin{tikzpicture}[scale=.72]
			\draw (0,0) circle (2cm);
			\draw [orange,line width=1] plot [smooth cycle] coordinates {(.3,0) (0,1.5) (.7,1) (1.5,1) (1.3,.4)};
			\fill [orange!10] plot [smooth cycle] coordinates {(.3,0) (0,1.5) (.7,1) (1.5,1) (1.3,.4)};
			\draw [dashed] (0,0) ellipse (2cm and 0.5cm);
			\draw (2,0) arc[x radius=2cm, y radius=.5cm, start angle=0, end angle=-180];
			\draw (-2,1.5) node {$\bS^{n-1}$};
			\draw [<-] (-3,0) to[bend right] (-5,0);
			\draw (-4,.5) node {$g(K,\cdot)$};
			\draw [rouge,->] (0,0) -- (.4,.7);
			\draw [rouge,->] (0,0) -- (.2,1.2);
			\draw [rouge,->] (0,0) -- (1,.7);
			\draw [rouge,->] (0,0) -- (.6,.2);
			\draw [orange] (-.3,1) node {$U$};
		\end{tikzpicture}
	\end{center}
\end{defn}

\begin{thm}[Minkowski's theorem]
	Let $\varphi$ be a finite measure on $\sB(\bS^{n-1})$ such that 
	$$\int_{\bS^{n-1}}ud\varphi(u)=0$$ and $\varphi(s)<\varphi(\bS^{n-1})$ for each great subsphere.
	
	Then there exists a convex body $K$ unique (up to translation) such that $$S_{n-1}(K,\cdot)=\varphi\ .$$
\end{thm}

\begin{prr}
	See \cite{ref6} Section 7 for example.
\end{prr}

We can now define the Blaschke sum of two convex bodies.

\begin{defn}\label{bs}
	Consider two convex bodies $K$ and $L$. The Blaschke sum of $K$ and $L$ is the unique \footnote{up to translation} convex body $M$ such that 	
		$$S_{n-1}(M,\cdot)=S_{n-1}(K,\cdot)+S_{n-1}(L,\cdot)\ .$$
\end{defn}

\begin{rem}
	If $K=P, L=Q$ are polytopes then we denote by $u_1,...,u_N\in\bS^{n-1}$ a family of pairwise distinct vectors linearly spanning $\bR^n$ such that the exterior normal vector of any face of $P$ or $Q$ is an element of $\{u_1,...,u_N\}$. Moreover we let $f_1, ..., f_N$ and $g_1,...,g_N$ (we may have $f_i=0$ or $g_i=0$) be the positive real numbers defined by 
	$$f_i=\vo(F(P,u_i))\quad\text{and}\quad g_i=\vo(F(Q,u_i))$$
	where $F(P,u)$ is the face of $P$ with exterior normal vector $u$.
	
	Then there exists a unique polytope $R$ whose exterior normal vectors are contained in $\{u_1,...,u_N\}$ such that the volume of the face normal to $u_i$ is 
	$$\vo(F(R,u_i))=f_i+g_i\ .$$
	In particular the Minkowski theorem implies that the Blaschke sum of two polytopes $P$ and $Q$ is a polytope.
\end{rem}

\begin{ex2}
	The Blaschke sum of a $n$ dimensional cube of side-length $d$ with itself is  
		$$[0,d]^n\#[0,d]^n=[0,\sqrt[n-1]{2}d]^n\ .$$
\end{ex2}

\begin{prop}
	If $K$ and $L$ are $2$ dimensional-polytopes then their Blaschke sum and their Minkowski sum coincide, $K\#L=K+L$.
\end{prop}

\begin{prr}
	A point in an edge of the Minkowski sum $K+L$ is a sum of two points respectively in the edge with same exterior normal vector in $K$ and in $L$. So the edges of $K+L$ are of length the sum of the lengths of the edge with same exterior normal vector in $K$ and in $L$. This is the property defining uniquely the Blaschke sum $K\#L$.
\end{prr}

\subsection{Motivation}\label{motiv}

In this subsection we discuss our motivation for studying the Blaschke sum. The first one was highlighted by Lehmann and Xiao in \cite{ref10} Section 7.A.

\subsubsection{Volume in equalities}

It is proved in \cite{ref2} that if $K$ and $L$ are any convex bodies then we have 
$$\vo(K\#L)^{\frac{n-1}n}\geq\vo(K)^{\frac{n-1}n}+\vo(L)^{\frac{n-1}n}$$
with equality if and only if $K$ and $L$ are homothetic. This is called the log concavity of the (Euclidean) volume.

\begin{prop}\label{convincing}
	The volume function $\fM$ on curves is $\frac{n-1}n$-log concave on movable curves with $\fM>0$.
\end{prop}

\begin{prr}
	Consider two movable curve classes $\alpha_1$ and $\alpha_2$ with $\fM(\alpha_1),\fM(\alpha_2)>0$ then
	\begin{align*}	
		\fM(\alpha_1+\alpha_2)^{\frac{n-1}n}
		&=\inf_{\substack{A\text{ big and movable}\\\text{divisor class}}}\frac{A\cdot(\alpha_1+\alpha_2)}{\vo(A)^\frac1n}\\
		&\geq\inf_{\substack{A\text{ big and movable}\\\text{divisor class}}}\frac{A\cdot\alpha_1}{\vo(A)^\frac1n}+\inf_{\substack{A\text{ big and movable}\\\text{divisor class}}}\frac{A\cdot\alpha_2}{\vo(A)^\frac1n}\\ 
		&=\fM(\alpha_1)^{\frac{n-1}n}+\fM(\alpha_2)^{\frac{n-1}n}.
	\end{align*} \vspace{-16pt}
\end{prr}

\subsubsection{Surfaces}

\begin{prop}
	The inclusion 	
		\begin{equation}
			\tag{\ref{incmov}}\Delta_{Y_\bullet}(\alpha_1)\#\Delta_{Y_\bullet}(\alpha_2)\subseteq\Delta_{Y_\bullet}(\alpha_1+\alpha_2)
		\end{equation}
	holds for any projective surface. 
\end{prop}

\begin{prr}
	If $X$ is a surface, the Newton--Okounkov bodies for curves and divisors coincide and moreover the Blaschke sum coincides with the Minkowski sum. 
	
	Indeed \cite{ref17} gives the equality of the Minkowski sum and the Blaschke sum on polyhedra and Küronya, Lozovanu and Maclean proved that all Newton--Okounkov bodies in dimension $2$ are polyhedral (see Theorem B \cite{ref15}). 
	
	As a consequence, the inclusion $\Delta_{Y_\bullet}(\alpha_1)\#\Delta_{Y_\bullet}(\alpha_2)\subseteq\Delta_{Y_\bullet}(\alpha_1+\alpha_2)$ can be rewritten as 	
		$$\Delta_{Y_\bullet}(B_1)+\Delta_{Y_\bullet}(B_2)\subseteq\Delta_{Y_\bullet}(B_1+B_2)_ ,$$ 
	which holds by the following.
	
	If $f$ is a section of $\sO_X(lB_1)$ and $g$ is a section of $\sO_X(mB_2)$ then $f^mg^l$ is a section of $\sO_X\big(lm(B_1+B_2)\big)$, and 	
		$$v_{Y_\bullet}(B_1+B_2)(f^mg^l)=mv_{Y_\bullet}(B_1)(f)+lv_{Y_\bullet}(B_2)(g)$$ 
	where $v_{Y_\bullet}(D)$ are the valuations used in the construction of the Newton--Okounkov bodies $\Delta_{Y_\bullet}(D)$. We then have that 	
		$$\Gamma_{Y_\bullet}(B_1)+\Gamma_{Y_\bullet}(B_2)\subseteq\Gamma_{Y_\bullet}(B_1+B_2)\ ,$$ 
	and consequently 
		$$\Delta_{Y_\bullet}(B_1)+\Delta_{Y_\bullet}(B_2)\subseteq\Delta_{Y_\bullet}(B_1+B_2)\ .\vspace{-17pt}$$
\end{prr}

\subsubsection{Homothetic curve classes}

\begin{prop}\label{multiple}
	Assume that there exists $A=a^{r-1}\in\bR^*$ such that $\alpha_2=A\alpha_1$. Then 
		$$\Delta(\alpha_1)\#\Delta(\alpha_2)=\Delta(\alpha_1+\alpha_2)=(1+A)^{\frac 1{n-1}}\Delta(\alpha_1)\ .$$
\end{prop}

\begin{prr}
	If $\alpha_1=\langle L^{n-1}\rangle$ and $\alpha_2=\langle (aL)^{n-1}\rangle$ then $\Delta(\alpha_2)=a\Delta(\alpha_1)$. Moreover, we have \begin{center}
		$\alpha_1+\alpha_2=\langle \big((1+A)^{\frac1{n-1}}L\big)^{n-1}\rangle\quad\text{and}\quad\Delta(\alpha_1+\alpha_2)=(1+A)^{\frac 1{n-1}}\Delta(\alpha_1).$
	\end{center}
	If $u$ is a normal vector of a face of $\Delta(\alpha_1)$ then the volume of the face of $\Delta(\alpha_1+\alpha_2)$ with normal vector $u$ is exactly the sum of the volume of the faces of $\Delta(\alpha_1)$ and $\Delta(\alpha_2)$ with normal vector $u$ 
		$$\vspace{-24pt}	\vo\big(\Delta(\alpha_1+\alpha_2),u\big)=(1+A)\times\vo\big(\Delta(\alpha_1),u\big)=\vo\big(\Delta(\alpha_1),u\big)+\vo\big(\Delta(\alpha_2),u\big)\ .$$\vspace{6pt}	
\end{prr}

\begin{cor}
	Let $X$ be any projective variety with $N^1(X)_\bR$ generated by only one element. The inclusion (\ref{incmov}) holds automatically for every movable curve on $X$ positive with respect to $\fM$ (and is moreover an equality).
\end{cor}

\subsubsection{Toric varieties}

Toric varieties are known to offer computable examples. Consider a normal projective toric variety $X$ with torus $T$ and its associated fan $\Sigma$ (see \cite{ref20} for more information about toric varieties). We consider a compatible flag 
$$Y_\bullet:\quad X\supseteq D_1\supseteq D_1\cap D_2\supseteq...\supseteq D_1\cap...\cap D_d=\{pt\}$$ where $D_1,...,D_d$ are prime $T$- invariant divisors. 

In \cite{ref4} Proposition 6.1, Lazarsfeld and Musta\c{t}\u{a} show that the Newton--Okounkov body of any big $T$-invariant divisor with respect to the flag $Y_\bullet$ can be written as 
$$\Delta_{Y_\bullet}(L) = \Phi_{Y_\bullet,\bR}(\Delta(L))$$
where $\Delta(L)$ is the Newton polytope of $L$ and  $\Phi_{Y_\bullet,\bR}$ is the linear map obtained by tensorizing by $\bR$ the dual isomorphism $\Phi_{Y_\bullet}:M\to\bZ^d;\ u\mapsto (\langle u,u_i\rangle)_{1\leq i\leq d}$.

The same way, we define $\Delta(\alpha)$ to be the Newton polytope of $L_\alpha$ and we have 
$$\Delta_{Y_\bullet}(\alpha)=\Phi_{Y_\bullet,\bR}(\Delta(\alpha))\ .$$

\begin{prop}
	A movable curve class $\alpha$ on $X$ satisfies $\fM(\alpha)>0$ if and only if $\alpha$ is positive along a spanning set of
	rays of $\Sigma$.
\end{prop}

\begin{prr}
	See \cite{ref3} Lemma 4.1 and Theorem 4.2.
\end{prr}

\begin{prop}\label{toricm}
	If $\alpha_1,\alpha_2$ are movable curve classes on the toric variety $X$ such that $\fM(\alpha_i)>0$ for $i=1,2$ then we have 
		$$\Delta_{Y_\bullet}(\alpha_1)\#\Delta_{Y_\bullet}(\alpha_2)=\Delta_{Y_\bullet}(\alpha_1+\alpha_2)\ .$$ 
\end{prop}

\begin{prr}
	As $\Phi_{Y_\bullet}:M\to\bZ^d$ is an isomorphism, it is enough to prove that $$\Delta(\alpha_1)\#\Delta(\alpha_2)=\Delta(\alpha_1+\alpha_2)\ .$$ 
	This follows directly from Theorem $4.2$ of Lehmann and Xiao \cite{ref3}. Take $X$ to be a projective toric variety with invariant divisors $D_1,...,D_s$ corresponding to rays $\rho_1,...\rho_s$ in the fan $\Sigma$ generated by the vectors $u_1,...,u_s$.
	To a curve $\alpha$, we may associate by Minkowski's theorem a polytope $P_\alpha$ such that the volume of the face $F_i$ is 
		\begin{equation}\label{formula}
			f_i=\frac{(\alpha\cdot D_i)||u_i||}{(n-1)!}\ ,
		\end{equation}
	where $F_i$ the face of $P_\alpha$ with exterior normal vertor $u_i$.
	Lehmann and Xiao (\cite{ref3} Theorem 4.2) prove that if $\alpha$ is movable and $\fM(\alpha)>0$ then 	
		$$P_\alpha=\Delta(L_\alpha)$$
	which is by definition $\Delta(\alpha)$.
	Now consider two movable curve classes $\alpha_1,\alpha_2$ with $$\fM(\alpha_1),\fM(\alpha_2)>0\ .$$ Let $\alpha_3$ be their sum and let $f_i^1,f_i^2,f_i^3$ be the volume of the faces of $\Delta(\alpha_1),\Delta(\alpha_2),\Delta(\alpha_3)$ orthogonal to $u_i$ respectively. By (\ref{formula}) we have that
		$$f_i^3=f_i^1+f_i^2\quad\forall i\in\{1,2,3\}\ ,$$
	so that $\Delta(\alpha_3)$ is exactly $\Delta(\alpha_1)\#\Delta(\alpha_2)$.
\end{prr}

\section{Potential reduction of (\ref{incmov}) to complete intersection curve classes}\label{cont}

To prove (\ref{incmov}) in general it may be enough to prove it only for curve classes of the form 
$$\alpha=A^{n-1}\ ,$$ for $A$ an ample divisor class.
We prove that such a reduction is possible assuming the continuity of the Blaschke sum, that is to say the continuity of the maps
\begin{align*}
	f_\bQ:\{\text{convex bodies in }\bR^d\}&\to\{\text{convex bodies in }\bR^d\}\ ,\\
	P&\mapsto P\#Q
\end{align*}
for every convex body $Q$ in $\bR^d$ where the topology is the Hausdorff distance on $\bR^d$. 

This reduction requires us to consider Newton--Okounkov bodies constructed with respect to a general maximal rank valuation and not only valuations coming from a flag (see Definition 3.4 of \cite{ref26}). Let us first state our result.

\begin{prop}\label{reduction}
	Assume the continuity of the Blaschke sum. If for any projective complex variety $X$, any maximal rank valuation on $\bC(X)$ and any curves of the form $$\alpha_1=B_1^{n-1},\;\alpha_2=B_2^{n-1}$$ on $X$ with $B_1,B_2$ ample, the inclusion 
	\begin{equation}\tag{\ref{incmov}} 	
		\Delta_v(\alpha_1)\#\Delta_v(\alpha_2)\subseteq\Delta_v(\alpha_1+\alpha_2)
	\end{equation}
	holds, then the inclusion (\ref{incmov}) holds for any projective complex variety $X$, for any maximal rank valuation and for any movable curve classes.
\end{prop} 

\begin{rem}
	Boucksom constructed from any admissible flag satisfying $$Y_i|_{Y_{i-1}}\text{ is Cartier in }Y_{i-1}$$
	a valuation $v_{Y_\bullet}:\bC(X)\to\bZ^d$ called the flag valuation (see Example 2.17 of \cite{ref26}).
	
	In particular, it follows from Proposition \ref{reduction} that \begin{equation}\tag{\ref{incmov}} 	
		\Delta_{Y_\bullet}(\alpha_1)\#\Delta_{Y_\bullet}(\alpha_2)\subseteq\Delta_{Y_\bullet}(\alpha_1+\alpha_2)
	\end{equation} holds for any projective complex variety $X$, for any flag as above and for any movable curve classes.
\end{rem}

\begin{rem}\label{bou}
	To prove Proposition \ref{reduction}, we consider general movable curves and divisors $(\alpha_i,L_i)_{i\in\{1,2,3\}}$ such that
	$$\alpha_i=\langle L_i\rangle^{n-1}\ \text{ and }\ \alpha_3=\alpha_1+\alpha_2\ .$$ 
	In Definition 3.4 of \cite{ref26}, Boucksom generalized the notion of Newton--Okounkov bodies associated to a flag, to Newton--Okounkov bodies associated to any valuation on $\bC(X)$. These Newton--Okounkov bodies are defined up to translation.
	
	More precisely, fixing a line bundle $L$, a choice of section of $L$ provides an inclusion $H^0(X,L)\subset \bC(X)$ and thus a Newton--Okounkov body. A different choice of section leads to a translated Newton--Okounkov body.
	As $\bC(X)$ is a birational invariant, these Newton--Okounkov bodies, defined up to translation are fixed under birational maps.
	
	Consequently, we need to find some sufficiently good common Fujita approximations of $L_1,L_2,L_3$ such that the curves $\alpha_i$ are also approximated. This is the aim of the following lemma which is the argument of \cite{ref14} Theorem 6.22 and which we resume and adapt to the case of several divisors here for the convenience of the reader.
\end{rem}

\begin{lem}\label{fu}
	Consider movable Cartier divisors $L_1,L_2,L_3$ on a complex projective variety $X$. Fix an ample divisor $H$ on $X$. 
	
	Then, for any $m\in\bN_{>0}$, there exists a birational map
	$\pi_m : X_{\pi_m}\to X$ and ample divisor classes $A_{1,m}$, $A_{2,m}$ and $A_{3,m}$ such that the following properties hold.
	\begin{enumerate}
		\item[P1:]\label{en1} $\pi_m^*L_i- A_{i,m}$ is pseudo-effective for all $i\in\{1,2,3\}$\ ;
		\item[P2:]\label{en2} $\vo(A_{i,m})>\vo(L_i)-\frac1m$ for all $i\in\{1,2,3\}$\ ;
		\item[P3:]\label{en4} For any $\varepsilon>0$, for $m$ large enough and for all $i\in\{1,2,3\}$, we have $$\pi_m^*(\langle L_i^{n-1}\rangle-\varepsilon H^{n-1})\leq A_{i,m}^{n-1}\leq\pi_m^*((\langle L_i^{n-1}\rangle)\ .$$
	\end{enumerate}
\end{lem}

\begin{pro}
	Applying Proposition 3.7 of \cite{ref30} to each $L_i$, there exist effective divisors $G_i$ so that for any sufficiently large $m$ there is a smooth birational model $$\phi_{i,m}:X_{\phi_{i,m}}\to X$$
	and a big and nef divisor $N_{i,m}$ on $X_{\phi_{i,m}}$ such that
	\begin{equation}\tag{A}\label{eq}
		P_\sigma(\phi_{i,m}^*L_i)-\frac1m\phi_{i,m}^*G_i\leq N_{i,m}\leq P_\sigma(\phi^*_{i,m}L_i)\ ,
	\end{equation}
	where $P_\sigma(.)$ denotes the positive part of the $\sigma$-decomposition of divisors (see \cite{ref13}). 
	
	Because $G_i$ does not depend on $m$ in (\ref{eq}), we can further require that $G_i$ is effective and ample.
	Equation (\ref{eq}) implies that $(\phi_{i,m},N_{i,m})$ satisfies Property P2 and P1 respectively.
	
	Consider $\pi_m:X_{\pi_m}\to X$ the projection map of a common resolution of $\phi_{1,m}$, $\phi_{2,m}$ and $\phi_{3,m}$.
	Denote by $A_{i,m}$ a small perturbation of the pull back of $N_{i,m}$ to $X_{\pi_m}$ which is ample and still satisfies Properties P1 and P2.
	
	We now prove that it will also satisfy Property P3. Fix $\varepsilon>0$.
	Take $m$ large enough such that $L_i-\frac1mG_i$ are pseudo-effective.
	By Lemma 6.21 of \cite{ref14}, the positive product satisfies $$\pi^*\langle L^{n-1}\rangle=\langle\pi^*L^{n-1}\rangle\ ,$$
	for any birational morphism $\pi$ and big divisor $L$. Moreover by Proposition 4.13 of \cite{ref30}, the positive product is invariant under replacing a divisor by its positive part $$\langle L^{n-1}\rangle=\langle P_\sigma L^{n-1}\rangle\ .$$ We thus have that
	$$\phi_{i,m}^*\left\langle\left(L_i-\frac1mG_i\right)^{n-1}\right\rangle=\left\langle P_\sigma\left(\phi_{i,m}^*\left(L_i-\frac1mG_i\right)\right)^{n-1}\right\rangle\ ,$$
	and since $\phi_{i,m}^*G_i$ is big and nef, we have $P_\sigma(\phi_{i,m}^*L_i)-\frac1m\phi^*_{i,m}G_i\geq P_\sigma(\phi^*_{i,m}(L-\frac1mG_i))$ and consequently $$\phi_{i,m}^*\left\langle\left(L_i-\frac1mG_i\right)^{n-1}\right\rangle\leq \left\langle \left(P_\sigma(\phi_{i,m}^*L_i)-\frac1m\phi^*_{i,m}G_i\right)^{n-1}\right\rangle\leq N_{i,m}^{n-1}\ .$$
	Similarly, we obtain  $N_{i,m}^{n-1}\leq \langle P_\sigma(\phi_{i,m}^*L_i)^{n-1}\rangle=\phi_{i,m}^*\langle L_i^{n-1}\rangle.$
	
	Choosing $m$ sufficiently large, we may ensure that $\varepsilon H-(L_i^{n-1}-\langle(L-\frac1mG)^{n-1}\rangle)$ is movable
	by continuity of the positive product. Its pull back by $\phi_{i,m}$ is again movable and we have that \begin{equation}\tag{B}\label{eeq}
		\phi_{i,m}^{n-1}(\alpha_i-\varepsilon H)\leq N_{i,m}^{n-1}\leq\phi_{i,m}^*\alpha_i\ .
	\end{equation}
	To prove that (\ref{eeq}) remains true on any higher birational model $\pi$, it is enough to verify that (\ref{eq}) still holds under the pull back by $\pi$: the left hand side follows from $$\pi^*P_\sigma(L_i)\geq P_\sigma(\pi^*L)\ ,$$ and the right hand side comes from Proposition III.1.14 of \cite{ref13}. 
	
	By continuity of the positive product,  (\ref{eeq}) stays valid under small perturbation.
	We have thus proved that $(\pi_m,A_{1,m},A_{2,m},A_{3,m})$ satisfy P3.\qed
\end{pro}

\begin{prPr}[\ref{reduction}]
	We start with movable curve classes $\alpha_1$, $\alpha_2$ and $\alpha_3=\alpha_1+\alpha_2$ which are non negative with respect to $\fM$ and which induce by Theorem \ref{mzd} three big and movable divisor classes $L_1$, $L_2$ and $L_3$ satisfying 
	$$\alpha_i=\langle L_i^{n-1}\rangle\ .$$
	Fix an ample divisor $H$ on $X$ for which we apply Lemma \ref{fu}.
	
	By Property P1, considering the canonical section $s$ of $\pi_m^*L_i- A_{i,m}$, each section of $\sO_{X_{\pi_m}}(kA_{i,m})$ multiplied by $s$ gives rise to a section of $\sO_{X_{\pi_m}}(k\pi_m^*L_i)$ and we then have that 
	\begin{center}
		$\Delta_v(\pi_m^*L_i)\supseteq\Delta_v(A_{i,m})\quad$ for $i=1,2,3$ .
	\end{center}

	By Boucksom's construction (see Remark \ref{bou}), translation classes of Newton-Okounkov bodies are birational invariants. Thus the Newton--Okounkov bodies of $\pi_m^*L_i$ coincide for all $m$, up to translation. 
	
	Taking into account Property P2, Proposition \ref{c0} which indicates the continuity of Newton--Okounkov bodies, we obtain that $\Delta_v(A_{i,m})$ converge to the convex body $\Delta_v(\pi_m^*L_i)=\Delta_v(\alpha_i)$.\\
	Using Property P3, it follows that $A_{1,m}^{n-1}+A_{2,m}^{n-1}$ and $A_{3,m}^{n-1}$ both converge to $\alpha_3$ when $m$ tends to $\infty$.
	
	Finally we supposed that the inclusion $\Delta_v(A_{1,m}^{n-1}+A_{2,m}^{n-1})\supseteq\Delta_v(A_{1,m})\#\Delta_v(A_{2,m})$ holds, this induces the inclusion for the movable classes $\alpha_1,\alpha_2,\alpha_3$
	$$\Delta_v(\alpha_3)=\lim_{m\to \infty}\Delta_v(A_{1,m}^{n-1}+A_{2,m}^{n-1})\supseteq\lim_{m\to\infty}(\Delta_v(A_{1,m})\#\Delta_v(A_{2,m}))=\Delta_v(\alpha_1)\#\Delta_v(\alpha_2)\ .$$	
\end{prPr}

\section{Newton--Okounkov bodies for curve classes on projective bundles over curves}\label{nob,cones}

In this section, we plan to give a complete description of Newton--Okounkov bodies of curve classes on a projective bundle over a curve. First we will recall some generalities on projective bundles over curves and the cones of divisor and curve classes. 

In \cite{ref1} Montero calculates the precise form of Newton--Okounkov bodies of divisors on projective vector bundles over curves. We summarise the relevant facts and explain how to see the positivity of a divisor by means of its Newton--Okounkov body.

\subsection{Generalities on projective bundles over curves}

\begin{defn}[See \cite{ref12} Page 160]
	Given a vector bundle $E$ on a curve $C$, the projectivisation of $E$ is $X=
	\Proj(\Sym^\bullet E)$ where the symmetric algebra $\Sym^\bullet E$ is the graded $\sO_C$-algebra given by $$\Sym^\bullet E(U)=\bigoplus_{m\in\bN}H^0(U,E^{\otimes m})\ .$$
\end{defn}

The power of this definition is in the following proposition.

\begin{prop}\label{grothendieckdefn}
	Denote the natural bundle map by $\pi:X\to C$. The variety $X$ carries a natural line bundle $\sO_X(1)$ of quotients by hyperplanes in $F=\pi^{-1}(p)\subseteq X$, satisfying
		$$\pi_*(\sO_X(k))=\Sym^kE\quad\text{for all }k\in\bN\ .$$
	 
\end{prop}

Let us recall the intersection ring of a projective bundle over a curve.

\begin{prop}
Set $\chi=\sO_X(1)$ and let $f=\pi^{-1}(q)$ be the fiber of a point.
Every divisor on $X$ can be written in the form $D=a(\chi-tf)$ for some real $a$ and $t$. More precisely, the intersection ring of $X$ $$A(X)=\bigoplus_{i=0}^r A^i(X)$$ is a graded $\bR$-algebra with multiplication induced by the intersection form and generated in degree $1$ by $\chi$ and $f$ with relations 
$$f^2=0,\quad\chi^r=d\cdot [pt]\quad\text{and}\quad\chi^{r-1}\cdot f=[pt]\ ,$$ where $[pt]$ denotes the class of a point and $d$ is the degree of first Chern class $c_1(E)$. 
\end{prop}

\begin{rem}\label{splitA(X)}
If $X$ is a projectivized splitting vector bundle $$\bP(\sO_{\bP_1}(a_1)\oplus...\oplus\sO_{\bP_1}(a_r))$$ over $\bP^1$ then $X$ is a toric variety and we can recover its intersection ring using toric theory (see Example 7.3.5 and Theorem 12.5.3 of \cite{ref20}).
\end{rem}

A way to make the study of vector bundles easier is to only look at semistable vector bundles. This requires the definition of slope.

\begin{defn}
	Consider a curve $C$, a vector bundle $(E,\pi)$ on $C$ of rank $r$ and degree $d=\deg(c_1(E))$ and consider $X=\Proj(\Sym^\bullet E)$ its projectivisation. 
	
	We denote by $\mu(E)$ the slope $\mu(E)=d/r$ of $E$ (See \cite{ref12} Page 52).
\end{defn}

\begin{defn}
	A semistable vector bundle is a vector bundle $E$ such that for every subbundle $Y\subseteq E$ we have $$\mu(Y)\leq\mu(E)\ .$$
\end{defn}

The Harder--Narasimhan filtration enables us to decompose $E$ into semistable vector bundles.

\begin{defn}
	The Harder--Narasimhan filtration of $E$ is the unique increasing filtration of $E$ by sub-bundles
		$$HN_\bullet(E) : 0=E_l\subseteq E_{l-1}\subseteq...\subseteq E_1\subseteq E_0=E$$
	such that each of the quotients $E_{i-1}/E_i$ satisfies the following conditions:
	\begin{enumerate}
		\item Each quotient $E_{i-1}/E_i$ is a semistable vector bundle\ ;
		\item $\mu(E_{i-1}/E_i)<\mu(E_i/E_{i+1})$ for all $i\in\{1,...,l-1\}\ .$
	\end{enumerate}
\end{defn}

\begin{nota}
	We will denote by $\mu_i$ the slope and $r_i$ the rank of the quotient $E_{i-1}/E_i$. We define numbers $\sigma_1\geq\sigma_2\geq...\geq\sigma_r$ by 
		$$(\sigma_1,...,\sigma_r)=(\underbrace{\mu_l,...,\mu_l}_{r_l\text{ times}},\underbrace{\mu_{l-1},...,\mu_{l-1}}_{r_{l-1}\text{ times}},...,\underbrace{\mu_1,...,\mu_1}_{r_1\text{ times}})\ .$$
\end{nota}

\subsection{Cones of divisor and curve classes on projective bundles over curves}\label{allcones}

A precise description of the cones of cycle classes was given by Fulger in \cite{ref8} and Fulger--Lehmann in \cite{ref14}. We will concentrate on the cones of divisor and curve classes
	$$\textup{Nef}_k(X)\subset\textup{Mov}_k(X)\subset\textup{Eff}_k(X)$$ 
with $k\in\{1,n-1\}$.

The description of the nef cones is due to Fulger (\cite{ref8} Lemma 2.1).
The pseudo-effective cone of divisors was computed by Nakayama (\cite{ref13} Corollary IV.3.8.). Fulger generalised it in \cite{ref8} Theorem 1.1 computing the pseudo-effective cone for any codimension. 
In \cite{ref14} Proposition 7.1, Fulger and Lehmann gave a similar description of the movable cone of cycle classes of $X$. 

\begin{prop}\label{posi}
	The effective, movable and nef cones of divisor and curve classes on a projective bundle $\bP(E)$ over a curve $C$ as above are given by
	$$\left\{\begin{array}{rcl}
		\Eff^1(X)&=&\langle f\;,\;\;\chi-\sigma_1f\rangle\ ,\\
		\mo^1(X)&=&\langle f\;,\;\;\chi-\sigma_2f\rangle\ ,\\ \Nef^1(X)&=&\langle f\;,\;\;\chi-\sigma_rf\rangle\ ,\\ \Eff_1(X)&=&\langle \chi^{r-2}\cdot f\;,\;\;\chi^{r-1}-(d-\sigma_r)\chi^{r-2}\cdot f\rangle\ ,
		\\ \mo_1(X)&=&\langle \chi^{r-2}\cdot f\;,\;\;\chi^{r-1}-(d-\sigma_1)\chi^{r-2}\cdot f\rangle\ ,\\ \Nef_1(X)&=&\langle \chi^{r-2}\cdot f\;,\;\;\chi^{r-1}-(r-1)\sigma_r\chi^{r-2}\cdot f\rangle\ .
	\end{array}\right.$$
\end{prop}

These data are summarized in the coming pictures. Let us first introduce a last cone of positivity: $C(X)$.

\begin{defn}
	We define $C(X)$ to be the cone generated by complete intersections of a unique nef divisor 
	$$C(X)=\overline{\textup{Cone}(\langle B^{r-1}\;|\,B\text{ big and nef}\rangle)}\ .$$
\end{defn}

In the case of projective bundles over curves, this cone can be described explicitly.

\begin{prop}
	The cone $C(X)$ is the set of curves of the form $$\chi^{r-1}-s\chi^{r-2}\cdot f$$ with $s<(r-1)\sigma_r$. Moreover $C(X)$ coincides with the complete intersection cone.
	\begin{align*}
		C(X)&=\langle \chi^{r-1}-(r-1)\sigma_r\chi^{r-2}\cdot f,\chi^{r-2}\cdot f\rangle\\
		&=\overline{\textup{Cone}(\langle B_1\cdot...\cdot B_{r-1}\;|\,B_i\text{ big and nef for all }i\rangle)}
	\end{align*}
\end{prop}

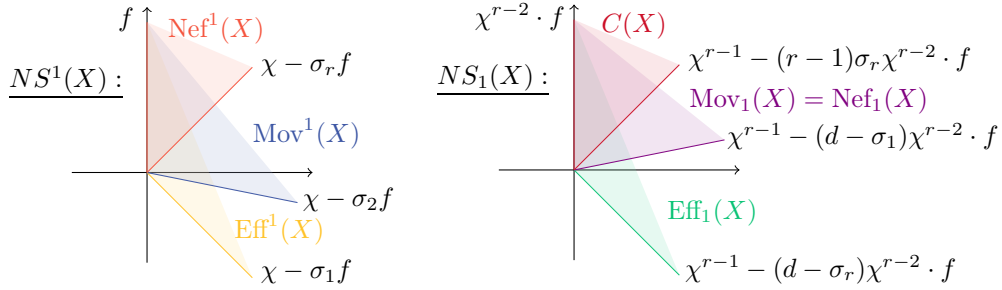
\begin{center}
	\begin{tikzpicture}[scale=0.4]
		\draw [color=noir,->] (-2.5,0) -> (5.5,0);
		\draw [color=noir,->] (0,-3) -> (0,5.5);
		\draw [color=jaune] (0,0) -- (3.5,-3.5); 
		\fill [color=jaune!50,opacity=0.2] (0,0) -- (0,5) -- (3.5,-3.5) -- cycle ;
		\draw [color=bleu] (0,0) -- (5,-1);
		\fill [color=bleu!50,opacity=0.2] (0,0) -- (0,5) -- (5,-1) -- cycle;
		\draw [color=orange] (0,0) -- (3.5,3.5);
		\draw [color=orange] (0,0) -- (0,5);
		\fill [color=orange!50,opacity=0.2] (0,0) -- (0,5) -- (3.5,3.5) -- cycle;
		\draw [color=jaune] (4.4,-1.9) node {$\Eff^1(X)$};
		\draw [color=bleu] (5.3,1.3) node {$\mo^1(X)$};
		\draw [color=orange] (2.3,4.8) node {$\Nef^1(X)$};
		\draw [color=noir] (-0.7,5.1) node {$f$};
		\draw [color=noir] (5.3,3.6) node {$\chi-\sigma_rf$};
		\draw [color=noir] (6.7,-0.9) node {$\chi-\sigma_2f$};
		\draw [color=noir] (5.3,-3.3) node {$\chi-\sigma_1f$};
		\draw [color=noir] (-2.7,3) node {$\uline{NS^1(X):}$};
	\end{tikzpicture}
	$\;$
	\begin{tikzpicture}[scale=0.4]
		\draw [color=noir,->] (-2.5,0) -> (5.5,0);
		\draw [color=noir,->] (0,-3) -> (0,5.5);
		\draw [color=vert] (0,0) -- (3.5,-3.5); 
		\fill [color=vert!50,opacity=0.2] (0,0) -- (0,5) -- (3.5,-3.5) -- cycle ;
		\draw [color=violet] (0,0) -- (5,1);
		\fill [color=violet!50,opacity=0.2] (0,0) -- (0,5) -- (5,1) -- cycle;
		\draw [color=rouge] (0,0) -- (3.5,3.5);
		\draw [color=rouge] (0,0) -- (0,5);
		\fill [color=rouge!50,opacity=0.2] (0,0) -- (0,5) -- (3.5,3.5) -- cycle;
		\draw [color=vert] (4.6,-1.4) node {$\Eff_1(X)$};
		\draw [color=violet] (7.8,2.3) node {$\mo_1(X)=\Nef_1(X)$};
		\draw [color=rouge] (2,4.8) node {$C(X)$};
		\draw [color=noir] (-1.8,5.1) node {$\chi^{r-2}\cdot f$};
		\draw [color=noir] (8.5,3.7) node {$\chi^{r-1}-(r-1)\sigma_r\chi^{r-2}\cdot f$};
		\draw [color=noir] (9.6,1.1) node {$\chi^{r-1}-(d-\sigma_1)\chi^{r-2}\cdot f$};
		\draw [color=noir] (8.2,-3.3) node {$\chi^{r-1}-(d-\sigma_r)\chi^{r-2}\cdot f$};
		\draw [color=noir] (-2.7,3) node {$\uline{NS_1(X):}$};
	\end{tikzpicture}
\end{center}

\begin{prr}
	Any nef divisor $B$ has the form $B=a(\chi-tf)$ where $t\leq\sigma_r$. It follows that	
	$$B^{r-1}=a^{r-1}\left(\chi^{r-1}-(r-1)t\chi^{r-2}\cdot f\right)$$ 
	and $(r-1)t\leq(r-1)\sigma_r$. Conversely if $s<(r-1)\sigma_r$ then we can write 
	$$\chi^{r-1}-s\chi^{r-2}\cdot f=B^{r-1}\ ,$$
	where $B=\chi-\frac{s}{r-1}f$ is a big and nef divisor. 
	
	For the second part of the proposition, each $B_i$ can be written in the form $$B_i=a_i(\chi-t_if)$$ with $t_i<\sigma_r$ and we have that
	$$B_1\cdot...\cdot B_{r-1}=\prod_{i=1}^{r-1}a_i\left(\chi^{r-1}-\sum_{i=1}^{r-1}t_i\chi^{r-2}\cdot f\right)$$
	where $\sum_{i=1}^{r-1}t_i<(r-1)\sigma_r$.
\end{prr}

\begin{rem}
	By \cite{ref16} Theorem 0.2, the movable and the nef cone of curve classes coincide
	$$\mo_1(X)=\Nef_1(X)\ .$$ 
\end{rem}

\begin{rem}
	Recalling that the intersection ring of $X$ is 
	$$A(X)=\bZ[\chi,f]/(f^2,\chi^{r+1},\chi^r-d\chi^{r-1}\cdot f)\ ,$$
	the form of the cones of curve classes follows from the form of the cones of divisor classes. Indeed by \cite{ref16} Theorem 0.2, the dual cone to the pseudo-effective cone of divisors $\Eff^1(X)=\langle f,\chi-\sigma_1f\rangle$ is the movable cone of curves and
	$$\left\{\begin{array}{rcl}
		f\cdot (\chi^{r-1}-t\chi^{r-2}\cdot f)&=& \chi^{r-1}\cdot f\\
		(\chi-\sigma_1f)\cdot (\chi^{r-1}-t\chi^{r-2}\cdot f)&=&(d-t+\sigma_1)\chi^{r-2}\cdot f	\end{array}\right.\ .$$
	In the same way we may deduce from the duality between the effective cone of curves and the cone of nef divisors the form of the latter.
\end{rem}

\subsection{Newton--Okounkov bodies of divisor classes}

In \cite{ref1}, Montero give the form of the Newton--Okounkov body of any divisor associated to a linear flag on a projective bundle over a curve. We recall everything here.

\begin{defn}
	A complete flag of subvarieties $Y_\bullet$ on the
	projective vector bundle $\bP(E)$ is called a linear flag if for some point $q\in C$ and $p\in\bP(E)$
	\begin{center}
		$Y_0=\bP(E)\ ,\quad Y_1=\pi^{-1}(q)\simeq\bP^{r-1}\ ,$\\
		$Y_i\simeq\bP^{r-i} \text{ is a linear subspace of } Y_{i-1}\ \forall i\in\{1,...,r\}\quad\text{and}\quad Y_r=\{p\}\ .$
	\end{center}
\end{defn}

We need the notion of complete linear flag $Y_\bullet^{\textup{HN}}$ compatible with the filtrations of Harder--Narasimhan. Consider the Harder--Narasimhan filtration of $E$
$$HN_\bullet(E) : 0=E_l\subseteq E_{l-1}\subseteq...\subseteq E_1\subseteq E_0=E\ .$$
There is a (possibly partial) flag of linear subvarieties
$$\bP((E/E_1)|_q)\subseteq\bP((E/E_2)|_q)\subseteq...\subseteq\bP((E/E_{l-1})|_q)\subseteq\bP(E|_q)=\pi^{-1}(q)\subseteq\bP(E)\ .$$
We will consider linear flags that are compatible with the Harder--Narasimhan filtration of $E$ in the sense that they complete the previous flag.

\begin{defn}
	A linear flag $Y_\bullet$ on $\bP(E)$ over $q\in C$ is said to be compatible with the Harder--Narasimhan filtration of $E$ if
	$$Y_{\rk E_i+1} = \bP((E/E_i)|_q) \simeq \bP^{r-\rk E_{i-1}}\subseteq\bP(E)\quad\text{for every }i\in\{1,...,l\}\ .$$
\end{defn}

We may decompose the full flag variety parameterising all complete linear flags into Schubert cells. The form of a Newton--Okounkov body associated to a complete linear flag will depend on the Schubert cell of the flag. 

\begin{defn}
	If we denote by $\bF_r$ the full flag variety parameterising all complete linear flags on $\pi^{-1}(q)\overset{f}\simeq\bP^{r-1}$, then there is a decomposition of $\bF_r$ into Schubert cells
	$$\bF_r=\bigsqcup_{\omega\in\fS_r}\Omega_\omega$$
	defined as follows.
	
	Consider homogeneous coordinates $[x_1:...,:x_r]$ on $\bP^{r-1}$ and let $Y_\bullet^\omega$ be the complete linear flag defined by
	$$Y_i^\omega=f_*\{x_1=...=x_i=0\}\subset Y_1=\pi^{-1}(q)\ .$$
	There is an action of $PGL_r(\bC)$ on $\bF_r$ via the natural action on the standard basis points $e_1,...,e_r\in\bP^{r-1}$. 
	The Schubert cell $\Omega_\omega$ is defined to be the orbit 
	$$\Omega_\omega=B\cdot Y_\bullet^\omega\ ,$$
	where $B$ is the subgroup of $PGL_r(\bC)$ that fixes a reference flag $Y_\bullet^{\textup{HN}}$.
	
	We say that a complete linear flag $Y_\bullet$ on $\bP(E)$ over $q\in C$ belongs to a Schubert cell $\Omega_\omega$ if the induced complete linear flag $Y_\bullet|_{Y_1}$ belongs to $\Omega_\omega$.
\end{defn}

To be able to define the decomposition into Schubert cells we needed a reference flag $Y_\bullet^{\textup{HN}}$. The following theorem of Montero tells us that as long as the reference flag is compatible with the Harder--Narasimhan filtration, Newton--Okounkov bodies do not depend on the choice of reference flag.

\begin{thm}\label{thmdef}
	Let $X$ be the projectivisation of a vector bundle $E$ on a curve $C$ and $Y_\bullet$ be any linear flag. Consider also a reference flag compatible with the Harder--Narasimhan filtration of $E$ and the decomposition of the full flag variety on $Y_1$ into Schubert cells 
	$$\bF_r=\bigsqcup_{\omega\in\fS_r}\Omega_\omega\ .$$
	
	Then the Newton--Okounkov body of $X,D=\chi-tf,Y_\bullet$ is of the form
	$$\Delta_{Y_\bullet}(D)=\Big\{(\nu_1,...,\nu_r)\in[0,+\infty[\times\Delta_{r-1}\,\Big|\;\nu_1+\sum_{i=2}^r\nu_i(\sigma_{\omega(r)}-\sigma_{\omega(i-1)})\leq\sigma_{\omega(r)}-t\Big\}\ ,$$ 
	where $\Delta_{r-1}$ is the unitary simplex of dimension $r-1$ and for some permutation $\omega$ corresponding to the Schubert cell of the flag $Y_\bullet$.
\end{thm}

\begin{prr}
	See Theorem 5.8 and Corollary 5.9 \cite{ref1}.
\end{prr}

\begin{rem}\label{H}
	The hyperplane $$H:\nu_1+\sum_{i=2}^r\nu_i(\sigma_{\omega(r)}-\sigma_{\omega(i-1)})=\sigma_{\omega(r)}-t$$ 
	splits $\bR^r$ into two half spaces and the Newton--Okounkov body $P$ is $[0,+\infty[\times\Delta_{r-1}$ intersected with one of them. 
	
	Moreover, in Proposition B of \cite{ref1}, Montero notes that the vector bundle $E$ is semistable if and only if the $\sigma_i$'s are equal and $H$ is 'straight' and has equation $\nu_1=\sigma_{\omega(r)}-t.$ In particular, the Newton--Okounkov body associated to any effective divisor $D=\chi-tf$ is $$\Delta_{Y_\bullet}(D)=[0,\sigma_{\omega(r)}-t]\times\Delta_{r-1}\ .$$
\end{rem}

\begin{ex2}
	Here are representations of all Newton--Okounkov bodies for $E$ of low rank.\\
	\uline{$\rk(E)=2:$}
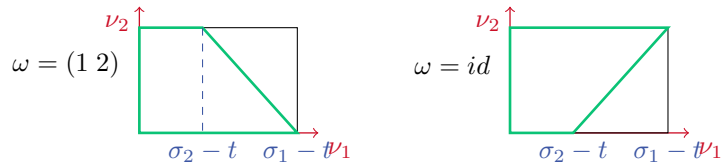
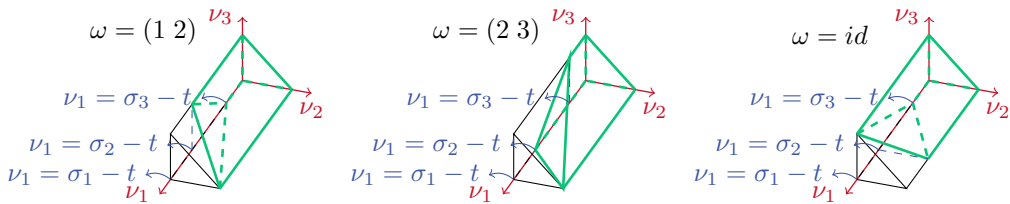
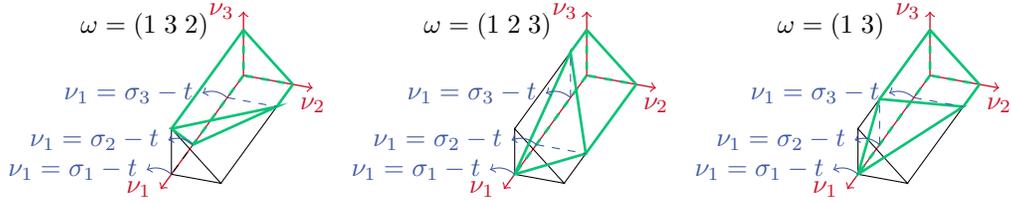
	\begin{center}
		\begin{tikzpicture}[scale=1.4]
			\draw[rouge,line width=.5pt,->] (0,0) -- (1.7,0) node[anchor=north west]{$\nu_1$};
			\draw[rouge,line width=.5pt,->] (0,0) -- (0,1.2) node[anchor=north east]{$\nu_2$};
			\draw (0,0) -- (1.5,0) -- (1.5,1) -- (0,1) -- cycle;
			\draw [vert, line width=1pt] (0,0) -- (1.5,0) -- (0.6,1) -- (0,1) -- cycle;
			\draw [color=bleu,dashed] (0.6,1) -- (0.6,0) node[below]{$\sigma_2-t$};
			\draw [color=bleu] (1.5,0) node[below] {$\sigma_1-t$};
			\draw (-.1,.65) node[left] {$\omega=(1\;2)$};
		\end{tikzpicture}
		$\quad$
		\begin{tikzpicture}[scale=1.4]
			\draw[rouge,line width=.5pt,->] (0,0) -- (1.7,0) node[anchor=north west]{$\nu_1$};
			\draw[rouge,line width=.5pt,->] (0,0) -- (0,1.2) node[anchor=north east]{$\nu_2$};
			\draw (0,0) -- (1.5,0) -- (1.5,1) -- (0,1) -- cycle;
			\draw [vert, line width=1pt] (0,0) -- (0.6,0) -- (1.5,1) -- (0,1) -- cycle;
			\draw [color=bleu] (0.6,0) node[below]{$\sigma_2-t$};
			\draw [color=bleu] (1.5,0) node[below] {$\sigma_1-t$}; 
			\draw (-.1,.65) node[left] {$\omega=id$};
		\end{tikzpicture}
	\end{center}
	\uline{$\rk(E)=3:$}
	\begin{center}
		\tdplotsetmaincoords{60}{110}
		\begin{tikzpicture}[tdplot_main_coords,scale=.7]
			\draw[rouge,line width=.5pt,->] (0,0,0) -- (4.6,0,0) node[left]{$\nu_1$};
			\draw[rouge,line width=.5pt,->] (0,0,0) -- (0,1.4,0) node[anchor=north]{$\nu_2$};
			\draw[rouge,line width=.5pt,->] (0,0,0) -- (0,0,1.4) node[anchor=east]{$\nu_3$};
			\draw[black,dashed] (0,0,0) -- (4,0,0);
			\draw[black,dashed] (0,0,1) -- (0,0,0) -- (0,1,0);
			\draw[black] (0,1,0) -- (0,0,1);
			\draw[black] (0,0,1) -- (4,0,1);
			\draw[black] (0,1,0) -- (4,1,0);
			\draw[black] (4,0,0) -- (4,0,1) -- (4,1,0) -- cycle;
			\draw[vert,dashed,line width=1pt] (0,0,0) -- (.9,0,0);
			\draw[vert,dashed,line width=1pt] (0,0,1) -- (0,0,0) -- (0,1,0);
			\draw[vert,line width=1pt] (0,1,0) -- (0,0,1);
			\draw[vert,line width=1pt] (0,0,1) -- (2.8,0,1);
			\draw[vert,line width=1pt] (0,1,0) -- (4,1,0);
			\draw[vert,dashed,line width=1pt] (4,1,0) -- (.9,0,0) -- (2.8,0,1) ;
			\draw[vert,line width=1pt] (2.8,0,1) -- (4,1,0);
			\draw [color=bleu,->] (4,0,0) to[bend right] (4,-.5,0) node [left] {$\nu_1=\sigma_1-t$};
			\draw [color=bleu,->] (2.8,0,0) to[bend right] (2.8,-.5,0) node [left] {$\nu_1=\sigma_2-t$};
			\draw [color=bleu,->] (.9,0,0) to[bend right] (.9,-.5,0) node [left] {$\nu_1=\sigma_3-t$};
			\draw [color=bleu,dashed] (2.8,0,0) -- (2.8,0,1);
			\draw (0,-2,.2) node[above] {$\omega=(1\;2)$};
		\end{tikzpicture}
		\begin{tikzpicture}[tdplot_main_coords,scale=.7]
			\draw[rouge,line width=.5pt,->] (0,0,0) -- (4.6,0,0) node[left]{$\nu_1$};
			\draw[rouge,line width=.5pt,->] (0,0,0) -- (0,1.4,0) node[anchor=north]{$\nu_2$};
			\draw[rouge,line width=.5pt,->] (0,0,0) -- (0,0,1.4) node[anchor=east]{$\nu_3$};
			\draw[black,dashed] (0,0,0) -- (4,0,0);
			\draw[black,dashed] (0,0,1) -- (0,0,0) -- (0,1,0);
			\draw[black] (0,1,0) -- (0,0,1);
			\draw[black] (0,0,1) -- (4,0,1);
			\draw[black] (0,1,0) -- (4,1,0);
			\draw[black] (4,0,0) -- (4,0,1) -- (4,1,0) -- cycle;
			\draw[vert,dashed,line width=1pt] (0,0,0) -- (2.8,0,0);
			\draw[vert,dashed,line width=1pt] (0,0,1) -- (0,0,0) -- (0,1,0);
			\draw[vert,line width=1pt] (0,1,0) -- (0,0,1);
			\draw[vert,line width=1pt] (0,0,1) -- (.9,0,1);
			\draw[vert,line width=1pt] (0,1,0) -- (4,1,0);
			\draw[vert,line width=1pt] (2.8,0,0) -- (.9,0,1) -- (4,1,0) -- cycle;
			\draw [color=bleu,->] (4,0,0) to[bend right] (4,-.5,0) node [left] {$\nu_1=\sigma_1-t$};
			\draw [color=bleu,->] (2.8,0,0) to[bend right] (2.8,-.5,0) node [left] {$\nu_1=\sigma_2-t$};
			\draw [color=bleu,->] (.9,0,0) to[bend right] (.9,-.5,0) node [left] {$\nu_1=\sigma_3-t$};		
			\draw [color=bleu,dashed] (.9,0,0) -- (.9,0,1);
			\draw (0,-2,.2) node[above] {$\omega=(2\;3)$};
		\end{tikzpicture}
		\begin{tikzpicture}[tdplot_main_coords,scale=.7]
			\draw[rouge,line width=.5pt,->] (0,0,0) -- (4.6,0,0) node[left]{$\nu_1$};
			\draw[rouge,line width=.5pt,->] (0,0,0) -- (0,1.4,0) node[anchor=north]{$\nu_2$};
			\draw[rouge,line width=.5pt,->] (0,0,0) -- (0,0,1.4) node[anchor=east]{$\nu_3$};
			\draw[black,dashed] (0,0,0) -- (4,0,0);
			\draw[black,dashed] (0,0,1) -- (0,0,0) -- (0,1,0);
			\draw[black] (0,1,0) -- (0,0,1);
			\draw[black] (0,0,1) -- (4,0,1);
			\draw[black] (0,1,0) -- (4,1,0);
			\draw[black] (4,0,0) -- (4,0,1) -- (4,1,0) -- cycle;
			\draw[vert,dashed,line width=1pt] (0,0,0) -- (.9,0,0);
			\draw[vert,dashed,line width=1pt] (0,0,1) -- (0,0,0) -- (0,1,0);
			\draw[vert,line width=1pt] (0,1,0) -- (0,0,1);
			\draw[vert,line width=1pt] (0,0,1) -- (4,0,1);
			\draw[vert,line width=1pt] (0,1,0) -- (2.8,1,0);
			\draw[vert,line width=1pt,dashed] (4,0,1) -- (.9,0,0) -- (2.8,1,0);
			\draw[vert,line width=1pt] (4,0,1) -- (2.8,1,0);
			\draw [color=bleu,->] (4,0,0) to[bend right] (4,-.5,0) node [left] {$\nu_1=\sigma_1-t$};
			\draw [color=bleu,->] (2.8,0,0) to[bend right] (2.8,-.5,0) node [left] {$\nu_1=\sigma_2-t$};
			\draw [color=bleu,->] (.9,0,0) to[bend right] (.9,-.5,0) node [left] {$\nu_1=\sigma_3-t$};		
			\draw [color=bleu,dashed] (2.8,0,0) -- (2.8,1,0);
			\draw (0,-2,.2) node[above] {$\omega=id$};
		\end{tikzpicture}
		\begin{tikzpicture}[tdplot_main_coords,scale=.7]
			\draw[rouge,line width=.5pt,->] (0,0,0) -- (4.6,0,0) node[left]{$\nu_1$};
			\draw[rouge,line width=.5pt,->] (0,0,0) -- (0,1.4,0) node[anchor=north]{$\nu_2$};
			\draw[rouge,line width=.5pt,->] (0,0,0) -- (0,0,1.4) node[anchor=east]{$\nu_3$};
			\draw[black,dashed] (0,0,0) -- (4,0,0);
			\draw[black,dashed] (0,0,1) -- (0,0,0) -- (0,1,0);
			\draw[black] (0,1,0) -- (0,0,1);
			\draw[black] (0,0,1) -- (4,0,1);
			\draw[black] (0,1,0) -- (4,1,0);
			\draw[black] (4,0,0) -- (4,0,1) -- (4,1,0) -- cycle;
			\draw[vert,dashed,line width=1pt] (0,0,0) -- (2.8,0,0);
			\draw[vert,dashed,line width=1pt] (0,0,1) -- (0,0,0) -- (0,1,0);
			\draw[vert,line width=1pt] (0,1,0) -- (0,0,1);
			\draw[vert,line width=1pt] (0,0,1) -- (4,0,1);
			\draw[vert,line width=1pt] (0,1,0) -- (.9,1,0);
			\draw[vert,line width=1pt] (2.8,0,0) -- (4,0,1) -- (.9,1,0) -- cycle;
			\draw [color=bleu,->] (4,0,0) to[bend right] (4,-.5,0) node [left] {$\nu_1=\sigma_1-t$};
			\draw [color=bleu,->] (2.8,0,0) to[bend right] (2.8,-.5,0) node [left] {$\nu_1=\sigma_2-t$};
			\draw [color=bleu,->] (.9,0,0) to[bend right] (.9,-.5,0) node [left] {$\nu_1=\sigma_3-t$};		
			\draw [color=bleu,dashed] (.9,0,0) -- (.9,1,0);
			\draw (0,-2,.2) node[above] {$\omega=(1\;3\;2)$};
		\end{tikzpicture}
		\begin{tikzpicture}[tdplot_main_coords,scale=.7]
			\draw[rouge,line width=.5pt,->] (0,0,0) -- (4.6,0,0) node[left]{$\nu_1$};
			\draw[rouge,line width=.5pt,->] (0,0,0) -- (0,1.4,0) node[anchor=north]{$\nu_2$};
			\draw[rouge,line width=.5pt,->] (0,0,0) -- (0,0,1.4) node[anchor=east]{$\nu_3$};
			\draw[black,dashed] (0,0,0) -- (4,0,0);
			\draw[black,dashed] (0,0,1) -- (0,0,0) -- (0,1,0);
			\draw[black] (0,1,0) -- (0,0,1);
			\draw[black] (0,0,1) -- (4,0,1);
			\draw[black] (0,1,0) -- (4,1,0);
			\draw[black] (4,0,0) -- (4,0,1) -- (4,1,0) -- cycle;
			\draw[vert,dashed,line width=1pt] (0,0,0) -- (4,0,0);
			\draw[vert,dashed,line width=1pt] (0,0,1) -- (0,0,0) -- (0,1,0);
			\draw[vert,line width=1pt] (0,1,0) -- (0,0,1);
			\draw[vert,line width=1pt] (0,0,1) -- (.9,0,1);
			\draw[vert,line width=1pt] (0,1,0) -- (2.8,1,0);
			\draw[vert,line width=1pt] (4,0,0) -- (.9,0,1) -- (2.8,1,0) -- cycle;
			\draw [color=bleu,->] (4,0,0) to[bend right] (4,-.5,0) node [left] {$\nu_1=\sigma_1-t$};
			\draw [color=bleu,->] (2.8,0,0) to[bend right] (2.8,-.5,0) node [left] {$\nu_1=\sigma_2-t$};
			\draw [color=bleu,->] (.9,0,0) to[bend right] (.9,-.5,0) node [left] {$\nu_1=\sigma_3-t$};		
			\draw [color=bleu,dashed] (.9,0,0) -- (.9,0,1);
			\draw [color=bleu,dashed] (2.8,0,0) -- (2.8,1,0);
			\draw (0,-2,.2) node[above] {$\omega=(1\;2\;3)$};
		\end{tikzpicture}
		\begin{tikzpicture}[tdplot_main_coords,scale=.7]
			\draw[rouge,line width=.5pt,->] (0,0,0) -- (4.6,0,0) node[left]{$\nu_1$};
			\draw[rouge,line width=.5pt,->] (0,0,0) -- (0,1.4,0) node[anchor=north]{$\nu_2$};
			\draw[rouge,line width=.5pt,->] (0,0,0) -- (0,0,1.4) node[anchor=east]{$\nu_3$};
			\draw[black,dashed] (0,0,0) -- (4,0,0);
			\draw[black,dashed] (0,0,1) -- (0,0,0) -- (0,1,0);
			\draw[black] (0,1,0) -- (0,0,1);
			\draw[black] (0,0,1) -- (4,0,1);
			\draw[black] (0,1,0) -- (4,1,0);
			\draw[black] (4,0,0) -- (4,0,1) -- (4,1,0) -- cycle;
			\draw[vert,dashed,line width=1pt] (0,0,0) -- (4,0,0);
			\draw[vert,dashed,line width=1pt] (0,0,1) -- (0,0,0) -- (0,1,0);
			\draw[vert,line width=1pt] (0,1,0) -- (0,0,1);
			\draw[vert,line width=1pt] (0,0,1) -- (2.8,0,1);
			\draw[vert,line width=1pt] (0,1,0) -- (.9,1,0);
			\draw[vert,line width=1pt] (4,0,0) -- (2.8,0,1) -- (.9,1,0) -- cycle;
			\draw [color=bleu,->] (4,0,0) to[bend right] (4,-.5,0) node [left] {$\nu_1=\sigma_1-t$};
			\draw [color=bleu,->] (2.8,0,0) to[bend right] (2.8,-.5,0) node [left] {$\nu_1=\sigma_2-t$};
			\draw [color=bleu,->] (.9,0,0) to[bend right] (.9,-.5,0) node [left] {$\nu_1=\sigma_3-t$};
			\draw [color=bleu,dashed] (.9,0,0) -- (.9,1,0);
			\draw [color=bleu,dashed] (2.8,0,0) -- (2.8,0,1);
			\draw (0,-2,.2) node[above] {$\omega=(1\;3)$};
		\end{tikzpicture}
	\end{center}
	From now on we will assume that our flag is in the Schubert cell corresponding to the permutation $\omega=(1\;2\;...\;r)$.
\end{ex2}

\begin{rem}
	Since $\Delta(aD)=a\Delta(D)$, the above enables us to calculate the form of any Newton--Okounkov body.
\end{rem}

It will be useful to think of Newton--Okounkov bodies as a succession of slices.

\begin{defn}\label{slice}
	The $i^{\text{th}}$ slice of the Newton--Okounkov body $\Delta(D)=\Delta\big(a(\chi-tf)\big)$ is the intersection $$S_i\Delta(D)=\left(\bR_+^{n-1}\times a[\sigma_{i+1}-t,\sigma_i-t]\right)\cap \Delta(D)\ .$$ 
	The final slice will be denoted $\fsl(\Delta(D))$.
\end{defn}	

\begin{ex2}
	Here are two examples in dimension $4$.
	\begin{center}		
		\tdplotsetmaincoords{60}{110}
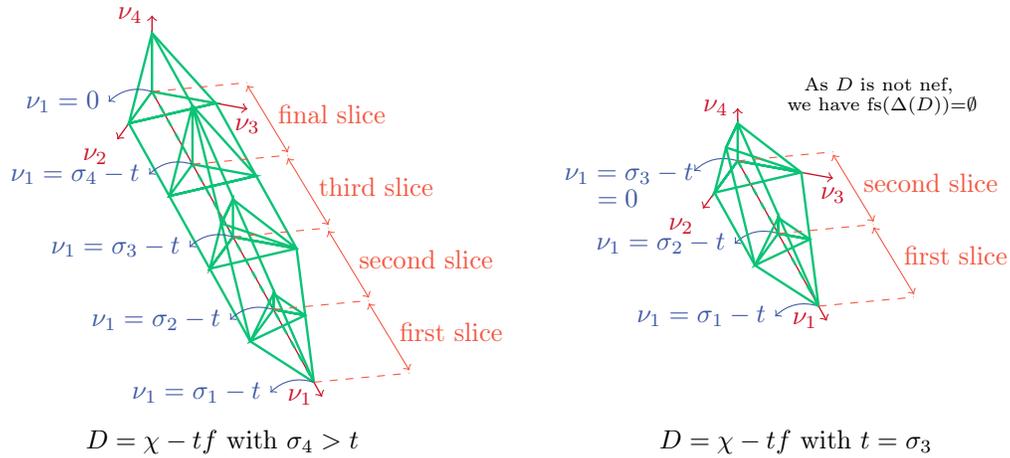
		\begin{tikzpicture}[tdplot_main_coords,scale=.9]
			\draw[rouge,line width=.5pt,->] (0,0,0) -- (4.2,4.2,-2.1) node[left]{$\nu_1$};
			\draw[rouge,line width=.5pt,->] (0,0,0) -- (1.5,0,0) node[anchor=north east]{$\nu_2$};
			\draw[rouge,line width=.5pt,->] (0,0,0) -- (0,1.5,0) node[anchor=north]{$\nu_3$};
			\draw[rouge,line width=.5pt,->] (0,0,0) -- (0,0,1.3) node[anchor=east]{$\nu_4$};
			\draw [vert,line width=.9pt,dashed] (0,0,0) -- (4,4,-2);
			\draw [vert,line width=.9pt] (4,4,-2) -- (4,3,-1.5) -- (1,0,0);
			\draw [vert,line width=.9pt] (4,4,-2) -- (2,3,-1) -- (0,1,0);
			\draw [vert,line width=.9pt] (4,4,-2) -- (1,1,.5) -- (0,0,1);
			\draw [vert,line width=.9pt] (0,0,0) -- (1,0,0) -- (0,1,0) -- (0,0,0) -- (0,0,1) -- (0,1,0) -- (1,0,0) -- (0,0,1);
			\draw [vert,line width=.9pt] (1,1,-.5) -- (2,1,-.5) -- (1,2,-.5) -- (1,1,-.5) -- (1,1,.5) -- (1,2,-.5) -- (2,1,-.5) -- (1,1,.5);
			\draw [vert,line width=.9pt] (2,2,-1) -- (3,2,-1) -- (2.5,2,-.5) -- (2,2,-.36) -- (2,3,-1) -- (3,2,-1) -- (2,2,-1) -- (2,3,-1) -- (2,2,-1) -- (2,2,-.36) -- (2.5,2,-.5) -- (2,3,-1);
			\draw [vert,line width=.9pt] (3,3,-1.5) -- (3,3,-1.23) -- (3,3.5,-1.5) -- (4,3,-1.5) -- (3,3,-1.5) -- (3,3.5,-1.5) -- (3,3,-1.5) -- (3,3,-1.23) -- (4,3,-1.5);
			\draw [vert,line width=.9pt] (4,3,-1.5) -- (1,1,.5);
			\draw [color=bleu,->] (0,0,0) to[bend right] (.5,-.5,0) node [left] {$\nu_1=0$};
			\draw [color=bleu,->] (1,1,-.5) to[bend right] (1.5,.5,-.5) node [left] {$\nu_1=\sigma_4-t$};
			\draw [color=bleu,->] (2,2,-1) to[bend right] (2.5,1.5,-1) node [left] {$\nu_1=\sigma_3-t$};
			\draw [color=bleu,->] (3,3,-1.5) to[bend right] (3.5,2.5,-1.5) node [left] {$\nu_1=\sigma_2-t$};
			\draw [color=bleu,->] (4,4,-2) to[bend right] (4.5,3.5,-2) node [left] {$\nu_1=\sigma_1-t$};
			\draw [color=orange,<->] (0,1.5,.4) -- (1,2.5,0);
			\draw [color=orange] (1,2.2,.6) node [right] {final slice};
			\draw [color=orange,<->] (1,2.5,-.1) -- (2,3.5,-.5);
			\draw [color=orange] (2,3.2,.1) node [right] {third slice};
			\draw [color=orange,<->] (2,3.5,-.6) -- (3,4.5,-1);
			\draw [color=orange] (3,4.2,-.4) node [right] {second slice};
			\draw [color=orange,<->] (3,4.5,-1.1) -- (4,5.5,-1.5);
			\draw [color=orange] (4,5.2,-.9) node [right] {first slice};
			\draw [orange,dashed] (0,0,0) -- (0,1.5,.45);
			\draw [orange,dashed] (1,1,-.5) -- (1,2.5,-.05);
			\draw [orange,dashed] (2,2,-1) -- (2,3.5,-.55);
			\draw [orange,dashed] (3,3,-1.5) -- (3,4.5,-1.05);
			\draw [orange,dashed] (4,4,-2) -- (4,5.5,-1.55);
		\end{tikzpicture}
		$\quad$
		\begin{tikzpicture}[tdplot_main_coords,scale=.9]
			\draw [color=bleu,->] (3,3,-1.5) to[bend right] (3.5,2.5,-1.5) node [left] {$\nu_1=\sigma_2-t$};
			\draw [color=bleu,->] (4,4,-2) to[bend right] (4.5,3.5,-2) node [left] {$\nu_1=\sigma_1-t$};
			\draw [color=bleu,->] (2,2,-1) to[bend right] (2.5,1.5,-1);
			\draw [color=bleu] (3.4,.8,-.7) node {$\nu_1=\sigma_3-t$};
			\draw [bleu] (2.3,.4,-1.8) node {$=0\;\;\;$};
			\draw[rouge,line width=.5pt,->] (2,2,-1) -- (4.2,4.2,-2.1) node[left]{$\nu_1$};
			\draw[rouge,line width=.5pt,->] (2,2,-1) -- (3.5,2,-1) node[anchor=north east]{$\nu_2$};
			\draw[rouge,line width=.5pt,->] (2,2,-1) -- (2,3.5,-1) node[anchor=north]{$\nu_3$};
			\draw[rouge,line width=.5pt,->] (2,2,-1) -- (2,2,-.1) node[anchor=east]{$\nu_4$};
			\draw [vert,line width=.9pt] (3,3,-1.5) -- (3,3,-1.23) -- (3,3.5,-1.5) -- (4,3,-1.5) -- (3,3,-1.5) -- (3,3.5,-1.5) -- (3,3,-1.5) -- (3,3,-1.23) -- (4,3,-1.5);
			\draw [vert,line width=.9pt] (2,2,-1) -- (3,2,-1) -- (2.5,2,-.5) -- (2,2,-.36) -- (2,3,-1) -- (3,2,-1) -- (2,2,-1) -- (2,3,-1) -- (2,2,-1) -- (2,2,-.36) -- (2.5,2,-.5) -- (2,3,-1);
			\draw [vert,line width=.9pt] (2,2,-.36) -- (4,4,-2);
			\draw [vert,line width=.9pt] (2,3,-1) -- (4,4,-2);
			\draw [vert,line width=.9pt] (3,2,-1) -- (4,3,-1.5) -- (4,4,-2);
			\draw [vert,line width=.9pt,dashed] (2,2,-1) -- (4,4,-2);
			\draw [vert, line width=.9pt] (2.5,2,-.5) -- (4,3,-1.5);
			\draw [color=orange,<->] (2,3.5,-.6) -- (3,4.5,-1);
			\draw [color=orange] (3,4.2,-.4) node [right] {second slice};
			\draw [color=orange,<->] (3,4.5,-1.1) -- (4,5.5,-1.5);
			\draw [color=orange] (4,5.2,-.9) node [right] {first slice};
			\draw [orange,dashed] (2,2,-1) -- (2,3.5,-.55);
			\draw [orange,dashed] (3,3,-1.5) -- (3,4.5,-1.05);
			\draw [orange,dashed] (4,4,-2) -- (4,5.5,-1.55);
			\draw [white] (0,0,-6.5) node {$g$};
			\draw (-1.5,3,-1.6) node {$\substack{\text{As $D$ is not nef, }\\ \text{we have }\fsl(\Delta(D))=\emptyset}$};
		\end{tikzpicture}
		\\	
		$D=\chi-tf$ with $\sigma_4>t$\hspace{4cm}$D=\chi-tf$ with $t=\sigma_3$
	\end{center} 
\end{ex2}

\begin{rem}
	The intersection of $\Delta(D)$ with the hyperplane $\nu_1=\sigma_j-t$ is the intersection of two simplexes 	
	$$\Big\{(\nu_2,...,\nu_r)\in\Delta_{r-1}\;\Big|\,\sum_{i=2}^r\nu_i\frac{\sigma_1-\sigma_i}{\sigma_1-\sigma_j}\leq1\Big\}\ .$$
\end{rem}

\begin{rem}\label{posit}
	We can translate the characterisations of the positivity (Proposition \ref{posi}) of a divisor in terms of its Newton--Okounkov body. 
	\begin{itemize}
		\item A divisor is nef if all the slices of its Newton--Okounkov body are non-empty. 
		\item A divisor is effective if its Newton--Okounkov body possesses a slice (i.e. is non-empty). 
		\item A divisor is movable if its Newton--Okounkov body contains at least the entire first slice. 
		In particular, if $\sigma_1=\sigma_2$ and $D=\chi-tf$ is an effective divisor then we have $\sigma_1=\sigma_2\geq t$, the first slice is then both empty and full and consequently $D$ is movable.
	\end{itemize}
		
	\tdplotsetmaincoords{60}{110}
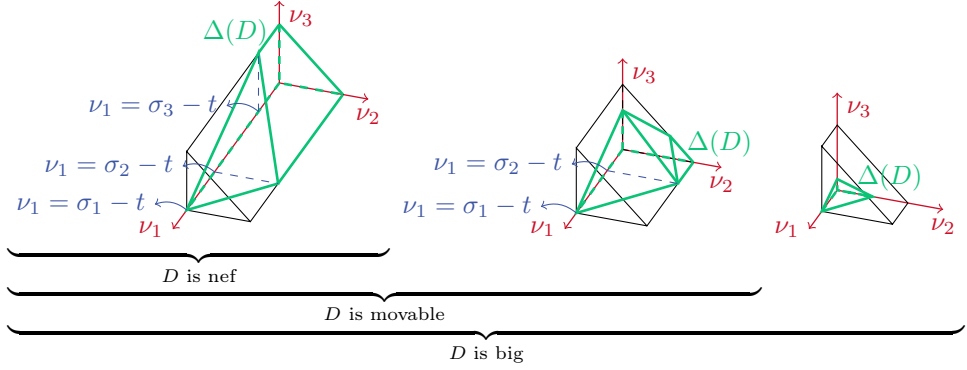
	\begin{center}$\underbrace{\underbrace{\underbrace{
					\begin{tikzpicture}[tdplot_main_coords,scale=.9]
						\draw[rouge,line width=.5pt,->] (0,0,0) -- (4.6,0,0) node[left]{$\nu_1$};
						\draw[rouge,line width=.5pt,->] (0,0,0) -- (0,1.4,0) node[anchor=north]{$\nu_2$};
						\draw[rouge,line width=.5pt,->] (0,0,0) -- (0,0,1.4) node[below right]{$\nu_3$};
						\draw (4,1,0) -- (4,0,0) -- (4,0,1) -- (4,1,0) -- (0,1,0) -- (0,0,1) -- (4,0,1);
						\draw[black,dashed] (0,0,0) -- (4,0,0);
						\draw[black,dashed] (0,0,1) -- (0,0,0) -- (0,1,0);
						\draw[vert,dashed,line width=1pt] (0,0,0) -- (4,0,0);
						\draw[vert,dashed,line width=1pt] (0,0,1) -- (0,0,0) -- (0,1,0);
						\draw[vert,line width=1pt] (0,1,0) -- (0,0,1);
						\draw[vert,line width=1pt] (0,0,1) -- (.9,0,1);
						\draw[vert,line width=1pt] (0,1,0) -- (2.8,1,0);
						\draw[vert,line width=1pt] (4,0,0) -- (.9,0,1) -- (2.8,1,0) -- cycle;
						\draw [color=bleu,->] (4,0,0) to[bend right] (4,-.5,0) node [left] {$\nu_1=\sigma_1-t$};
						\draw [color=bleu,->] (2.8,0,0) to[bend right] (2.8,-.5,0) node [left] {$\nu_1=\sigma_2-t$};
						\draw [color=bleu,->] (.9,0,0) to[bend right] (.9,-.5,0) node [left] {$\nu_1=\sigma_3-t$};		
						\draw [color=bleu,dashed] (.9,0,0) -- (.9,0,1);
						\draw [color=bleu,dashed] (2.8,0,0) -- (2.8,1,0); 
						\draw [color=vert] (-.3,-.8,.55) node {$\Delta(D)$};
				\end{tikzpicture}}_{D\text{ is nef}}
				\begin{tikzpicture}[tdplot_main_coords,scale=1]
					\draw[rouge,line width=.5pt,->] (0,0,0) -- (2.25,0,0) node[left]{$\nu_1$};
					\draw[rouge,line width=.5pt,->] (0,0,0) -- (0,1.4,0) node[anchor=north]{$\nu_2$};
					\draw[rouge,line width=.5pt,->] (0,0,0) -- (0,0,1.4) node[below right]{$\nu_3$};
					\draw (1.8,1,0) -- (1.8,0,0) -- (1.8,0,1) -- (1.8,1,0) -- (0,1,0) -- (0,0,1) -- (1.8,0,1);
					\draw[black,dashed] (0,0,0) -- (1.8,0,0);
					\draw[black,dashed] (0,0,1) -- (0,0,0) -- (0,1,0);
					\draw[vert,dashed,line width=1pt] (0,0,0) -- (1.8,0,0);
					\draw[vert,dashed,line width=1pt] (0,0,.6) -- (0,0,0) -- (0,1,0);
					\draw[vert,line width=1pt] (0,1,0) -- (0,.67,.34);
					\draw[vert,line width=1pt] (0,0,.6) -- (.6,1,0);
					\draw[vert,line width=1pt] (0,1,0) -- (.6,1,0);
					\draw[vert,line width=1pt] (1.8,0,0) -- (0,0,.6) -- (0,.67,.34) -- (.6,1,0) -- cycle;
					\draw [color=bleu,->] (1.8,0,0) to[bend right] (1.8,-.5,0) node [left] {$\nu_1=\sigma_1-t$};
					\draw [color=bleu,->] (.6,0,0) to[bend right] (.6,-.5,0) node [left] {$\nu_1=\sigma_2-t$};		
					\draw [color=bleu,dashed] (.6,0,0) -- (.6,1,0); 
					\draw [color=vert] (-.75,1.1,-.1) node {$\Delta(D)$};
			\end{tikzpicture}}_{D\text{ is movable}}
			\begin{tikzpicture}[tdplot_main_coords,scale=1]
				\draw[rouge,line width=.5pt,->] (0,0,0) -- (1.1,0,0) node[left]{$\nu_1$};
				\draw[rouge,line width=.5pt,->] (0,0,0) -- (0,1.5,0) node[anchor=north]{$\nu_2$};
				\draw[rouge,line width=.5pt,->] (0,0,0) -- (0,0,1.5) node[below right]{$\nu_3$};
				\draw (.6,1,0) -- (.6,0,0) -- (.6,0,1) -- (.6,1,0) -- (0,1,0) -- (0,0,1) -- (.6,0,1);
				\draw [vert, line width=1pt] (.6,0,0) -- (0,.5,0) -- (0,0,.17) -- cycle; 
				\draw [vert, line width=1pt,dashed] (.6,0,0) -- (0,0,0) -- (0,.5,0) -- (0,0,0) -- (0,0,.17);
				\draw [color=vert] (-.4,.6,.1) node {$\Delta(D)$};
		\end{tikzpicture}}_{D\text{ is big}}$
	\end{center}
\end{rem}

\subsection{Newton--Okounkov bodies of movable curve classes}

In this subsection we describe Newton--Okounkov bodies of movable $\fM$-positive curve classes. By Theorem \ref{mzd} every movable curve class $\alpha=\chi^{r-1}-s\chi^{r-2}\cdot f$ with $\fM(\alpha)>0$ can be uniquely written in the form 
	$$\alpha=\langle L^{r-1}\rangle,\quad\text{where $L$ is movable}\ .$$

\begin{prop}\label{nobm}
	The movable Zariski decomposition of a movable curve class $\alpha=\chi^{r-1}-s\chi^{r-2}\cdot f$ is 
		$$\alpha=\langle (\chi-tf)^{r-1}\rangle\ .$$
	where for some $t\leq\sigma_2$. Moreover, if $\alpha$ belongs to $C(X)$ then $t=\frac s{r-1}$.
	
	The Newton--Okounkov body of $\alpha$ is then 
		$$\Delta(\alpha)=\Delta(L)=\Delta\left(\chi-tf\right)\ .$$
\end{prop}

\section{The inclusion $\Delta(\alpha_1)\#\Delta(\alpha_2)\subseteq\Delta(\alpha_1+\alpha_2)$ on projective bundles over curves}\label{Inc}

Consider a projective bundle $X$ over a curve. In this section we would like to find conditions on curve classes $\alpha_1,\alpha_2$ under which the inclusion
	$$\Delta(\alpha_1)\#\Delta(\alpha_2)\subseteq\Delta(\alpha_1+\alpha_2)$$ 
holds using the movable Zariski decomposition.

We start by computing the Blaschke sum of the Okounkov bodies of two nef divisors.

\subsection{Blaschke sum of Newton--Okounkov bodies of nef divisors}

\begin{prop}\label{important!}
	Consider some Newton--Okounkov bodies $P=\Delta_{Y_\bullet}\big(\chi-t_1f\big)$ and $Q=\Delta_{Y_\bullet}\big(a(\chi-t_2f)\big)$ associated to big and nef divisors (i.e. with $\sigma_r\geq t_1,t_2$). 
	
	The Blaschke sum $R$ of $P$ and $Q$ is then given by 
	$$R=P\#Q=\Delta_{Y_\bullet}\big(b(\chi-t_3f)\big)\ ,$$ 
	where $t_3=\left(\frac{t_1+a^{r-1}t_2}{1+a^{r-1}}\right)$ and $b=(1+a^{r-1})^{\frac{1}{r-1}}$.
\end{prop}

Before starting the proof we define the common component and the final slice of a Newton--Okounkov body associated to a big and nef divisor.

\begin{nota}
	The volume of the face of $P$ with exterior normal vector $u$ will be denoted by $\vo(P,u)$. We denote by $$v=\begin{pmatrix}-1\\\sigma_1-\sigma_2\\\sigma_1-\sigma_3\\\vdots\\\sigma_1-\sigma_r\end{pmatrix}$$ the normal vector of the hyperplane $$H:\nu_1+\sum_{i=2}^r\nu_i(\sigma_{\omega(r)}-\sigma_{\omega(i-1)})=\sigma_{\omega(r)}-t$$ (see Remark \ref{H}).
		
	We denote by $P_i$ the $i^{\text{th}}$ slice of $P$.
\end{nota}

\begin{defn}
	Let $P=\Delta_{Y_\bullet}\big(a(\chi-tf)\big)$ be the Newton--Okounkov body associated to a big and nef divisor. We define the common component and the final slice of $P$ as $$\cc(P)=P\cap([\sigma_1-t,\sigma_r-t]\times\bR^{r-1})\quad\text{and}\quad\fsl(P)=P\cap([\sigma_r-t,0]\times\bR^{r-1})\ .$$
	
	With these definitions in mind, we may write $P$ as 	
		$$P=\textup{glueing}_{F(\cc(P),-\nu_1),F(\fsl(P),\nu_1)}\big(\cc(P);\fsl(P)\big)\ .$$ 
	In other words, the polytope $P$ is the union of $\cc(P)$ and $\tau(\fsl(P))$ where $\tau$ is the unique translation identifying $F(\cc(P),-\nu_1)$ and $F(\fsl(P),\nu_1)$.
	\tdplotsetmaincoords{60}{110}
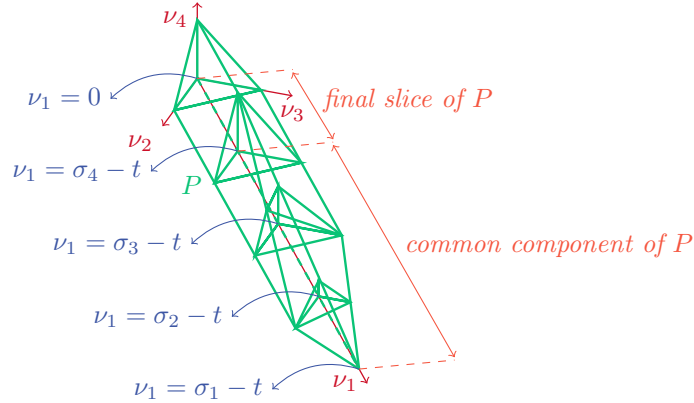
	\begin{center}\begin{tikzpicture}[tdplot_main_coords,scale=.9]
			\draw[rouge,line width=.5pt,->] (0,0,0) -- (4.2,4.2,-2.1) node[left]{$\nu_1$};
			\draw[rouge,line width=.5pt,->] (0,0,0) -- (1.5,0,0) node[anchor=north east]{$\nu_2$};
			\draw[rouge,line width=.5pt,->] (0,0,0) -- (0,1.5,0) node[anchor=north]{$\nu_3$};
			\draw[rouge,line width=.5pt,->] (0,0,0) -- (0,0,1.3) node[below left]{$\nu_4$};
			
			\draw [vert,line width=.9pt,dashed] (0,0,0) -- (4,4,-2);
			\draw [vert,line width=.9pt] (4,4,-2) -- (4,3,-1.5) -- (1,0,0);
			\draw [vert,line width=.9pt] (4,4,-2) -- (2,3,-1) -- (0,1,0);
			\draw [vert,line width=.9pt] (4,4,-2) -- (1,1,.5) -- (0,0,1);
			\draw [vert,line width=.9pt] (0,0,0) -- (1,0,0) -- (0,1,0) -- (0,0,0) -- (0,0,1) -- (0,1,0) -- (1,0,0) -- (0,0,1);
			\draw [vert,line width=.9pt] (1,1,-.5) -- (2,1,-.5) -- (1,2,-.5) -- (1,1,-.5) -- (1,1,.5) -- (1,2,-.5) -- (2,1,-.5) -- (1,1,.5);
			\draw [vert,line width=.9pt] (2,2,-1) -- (3,2,-1) -- (2.5,2,-.5) -- (2,2,-.36) -- (2,3,-1) -- (3,2,-1) -- (2,2,-1) -- (2,3,-1) -- (2,2,-1) -- (2,2,-.36) -- (2.5,2,-.5) -- (2,3,-1);
			\draw [vert,line width=.9pt] (3,3,-1.5) -- (3,3,-1.23) -- (3,3.5,-1.5) -- (4,3,-1.5) -- (3,3,-1.5) -- (3,3.5,-1.5) -- (3,3,-1.5) -- (3,3,-1.23) -- (4,3,-1.5);
			\draw [vert,line width=.9pt] (4,3,-1.5) -- (1,1,.5);
			\draw [color=bleu,->] (0,0,0) to[bend right] (1,-1,0) node [left] {$\nu_1=0$};
			\draw [color=bleu,->] (1,1,-.5) to[bend right] (2,0,-.5) node [left] {$\nu_1=\sigma_4-t$};
			\draw [color=bleu,->] (2,2,-1) to[bend right] (3,1,-1) node [left] {$\nu_1=\sigma_3-t$};
			\draw [color=bleu,->] (3,3,-1.5) to[bend right] (4,2,-1.5) node [left] {$\nu_1=\sigma_2-t$};
			\draw [color=bleu,->] (4,4,-2) to[bend right] (5,3,-2) node [left] {$\nu_1=\sigma_1-t$};
			\draw [color=orange,<->] (0,1.5,.4) -- (1,2.5,0);
			\draw [color=orange] (1,2.2,.6) node [right] {final slice of $P$};
			\draw [color=orange,<->] (1,2.5,-.1) -- (4,5.5,-1.5);
			\draw [color=orange] (3,4.2,-.4) node [right] {common component of $P$};
			\draw [vert] (3,1,0) node {$P$};
			\draw [orange,dashed] (0,0,0) -- (0,1.5,.45);
			\draw [orange,dashed] (1,1,-.5) -- (1,2.5,-.05);
			\draw [orange,dashed] (4,4,-2) -- (4,5.5,-1.55);
	\end{tikzpicture}\end{center}
\end{defn}

\begin{rem}
	Two nef divisors have the same common component up to homothety.
\end{rem}

\begin{prPr}[\ref{important!}]
	We claim that it is possible to decompose $P\#Q$ into the following gluing 
	$$P\#Q=\textup{glueing}\big(\cc(P)\#\cc(Q);\fsl(P)\#\fsl(Q)\big)\ .$$ 
	
	The main task is to prove that it is well-defined. As $$P=\Delta(\chi-t_1f)\quad\text{and}\quad Q=\Delta\big(a(\chi-t_2f)\big)$$ are nef,
	their common components are multiple of each other $\cc(Q)=a\times\cc(P)$. 
	
	Furthermore, by Proposition \ref{multiple} 	
		$$\cc(P)\#\cc(Q)=b\times\cc(P)\quad \text{with}\quad b=(1+a^{r-1})^{\frac{1}{r-1}}\ .$$ 
	
	The final slices are given by 	
		$$\fsl(P)=[0,\sigma_r-t_1]\times\Delta_{r-1}\:,\quad \fsl(Q)=a\times\big([0,\sigma_r-t_2]\times\Delta_{r-1}\big)$$
	where $\Delta_{r-1}$ is the $r-1$-dimensional simplex of vertices
		$$(0,0,...,0)\ ,\quad(1,0,...,0)\ ,\quad(0,1,0,...,0)\ ,\ \ldots\ ,\;(0,...,0,1)\ .$$ 
	
	By Minkowski's theorem the Blaschke product $\fsl(P)\#\fsl(Q)$ is the unique polytope $R$ whose faces have volume 	
		$$\begin{aligned}
			\vo\Big(F\big(R,u\big)\Big)&=\vo\Big(F\big(\fsl(P),u\big)\Big)+\vo\Big(F\big(\fsl(Q),u\big)\Big)\\
			&=\left\{
			\begin{aligned}
				\frac{1}{(r-2)!}+\frac{a^{r-1}}{(r-2)!}	&\quad\text{if } u=\pm\nu_1\ .\\
				\frac{\sigma_r-t_1}{(r-2)!}+\frac{a^{r-1}(\sigma_r-t_2)}{(r-2)!}	&\quad\text{if }u=-\nu_i\text{ for }i\in\{2,...,r\}\ .
			\end{aligned}\right.
		\end{aligned}$$
	
	Note that $\vo\Big(F\big(R,u\big)\Big)=\vo\Big(F\big(\fsl(P),u\big)\Big)+\vo\Big(F\big(\fsl(Q),u\big)\Big)$ for $u=\sum_{i=2}^r\nu_i$ follows from the previous equations: any polytope $S$ satisfies $\sum_{u\in\bS^1}F(S,u)u=0$.
	
	Consequently the Blaschke product $\fsl(P)\#\fsl(Q)$ is given by	
		$$ \fsl(P)\#\fsl(Q)=b\Big([0,\sigma_r-t_3]\times\Delta_{r-1}\Big)\quad
		\text{with }t_3=\frac{t_1+a^{r-1}t_2}{1+a^{r-1}}\;\;\text{and}\;\;b=(1+a^{r-1})^{\frac{1}{r-1}}\ .$$
	
	The faces $F(\cc(P)\#\cc(Q),-\nu_1)$ and $F(\fsl(P)\#\fsl(Q),\nu_1)$ coincide so that the polytopes $\cc(P)\#\cc(Q)$ and $\fsl(P)\#\fsl(Q)$ can be glued and 
		$$G:=\textup{glueing}\big(\cc(P)\#\cc(Q);\fsl(P)\#\fsl(Q)\big)$$ 
	is well defined. It is also convex. 
	
	It remains to prove that $G=P\#Q$. If $S$ is any polytope we denote by $S_j$ the $j^{\text{th}}$ slice of $S$. Let $F$ be a face of $G$ and let $u$ be the exterior normal vector of $F$. Then $u$ is either equal to 
	\begin{enumerate}
		\item $\nu_1$ and in that case \vspace{-1pt}	
		$$\vo(G,u)=\vo(G_1,u)=\vo(P_1,u)+\vo(Q_1,u)=\vo(P,u)+\vo(Q,u)\ ;$$
		\item $-\nu_1$ and in that case \vspace{-1pt}
		$$\vo(G,u)=\vo(G_r,u)=\vo(P_r,u)+\vo(Q_r,u)=\vo(P,u)+\vo(Q,u)\ ;$$
		\item $v$, $-\nu_2$, ..., $-\nu_r$ or $\sum_{j=2}^r\nu_j$ and in that case \vspace{-2pt}
		$$\vo(G,u)=\sum_{j=1}^r\vo(G_j,u)=\sum_{j=1}^r\vo(P_j,u)+\sum_{j=1}^r\vo(Q_j,u)=\vo(P,u)+\vo(Q,u)\ .$$
	\end{enumerate}
	It follows that $G=P\#Q$. Now since we have that
		$$\cc(P)\#\cc(Q)=b\times\cc(P)\quad\text{and}\quad\fsl(P)\#\fsl(Q)=b\Big([0,\sigma_r-t_3]\times\Delta_{r-1}\Big)\ ,$$
	the gluing $G$ is also equal to $\Delta\big(b(\chi-t_3f)\big).$ This completes the proof.
\end{prPr}

\subsection{Proof of the inclusion (\ref{incmov}) in the case $\alpha_1,\alpha_2\in C(X)$}

Every curve class $\alpha=c_1\chi^{r-1}+c_2\chi^{r-2}\cdot f$ can be written in the form
	$$\alpha=A^{r-1}\quad\text{where}\quad A=c_1^{\frac1{r-1}}\left(\chi-\frac{c_2}{(r-1)c_1}f\right)\ .$$
Let $A_1=\chi-t_1f$ and $A_2=a(\chi-t_2f)$ be two ample divisors and let $\alpha_1=A_1^{n-1}$ and $\alpha_2=A_2^{n-1}$ be their associated curve classes. 

The ample divisor $A_3$ satisfying $\alpha_3=A_3^{n-1}$ where $\alpha_3=\alpha_1+\alpha_2$ is equal to
	$$A_3=(1+a^{r-1})^{\frac1{r-1}}(\chi-t_3f)\quad\text{where}\quad t_3=\frac{t_1+a^{r-1}t_2}{1+a^{r-1}}\ .$$ 

We can now prove Theorem \ref{mthm}.

\begin{thm}\label{mthm}
	Let $X$ be a projective bundle over a curve. Consider two curves $\alpha_1,\alpha_2$ in $C(X)$ ie of the form $\alpha_i=A_i^{n-1}$ with $A_i$ ample divisor on $X$. 
	
	Then the inclusion (\ref{incmov}) holds and is an equality 
		$$\Delta(\alpha_1)\#\Delta(\alpha_2)=\Delta(\alpha_1+\alpha_2)\ .$$	
\end{thm}

\begin{prr}
	We may write $\alpha_1=a_1(\chi-t_1f)$ and $\alpha_2=a_2(\chi-t_2f)$. Set $a=\frac{a_2}{a_1}$ and apply Proposition \ref{important!} to $P=\Delta(\frac1{a_1}\alpha_1)=\Delta(\chi-t_1f)$ and $Q=\Delta(\frac1{a_1}\alpha_2)=\Delta\big(a(\chi-t_2f)\big)$. We obtain $$P\#Q=\Delta_{Y_\bullet}\big(b(\chi-t_3f)\big)\ ,$$ where $t_3=\frac{t_1+a^{r-1}t_2}{1+a^{r-1}}$ and $b=(1+a^{r-1})^{\frac{1}{r-1}}$.
	As $\alpha_1+\alpha_2=a_1(\frac1{a_1}\alpha_1+\frac1{a_1}\alpha_2)=ba_1(\chi-t_3f)$ we have 	
		$$\Delta(\alpha_1)\#\Delta(\alpha_2)=\Delta(\alpha_1+\alpha_2)\ .\vspace{-20pt}$$
\end{prr}

\newpage
\bibliographystyle{alpha} 
\bibliography{bib}

\Addresses

\end{document}